\documentclass[11pt]{article}
\usepackage[]{amsmath,amssymb}
\usepackage{cite}
\newtheorem{theorem}{Theorem}[section]
\newtheorem{lemma}[theorem]{Lemma}

\newtheorem{corollary}[theorem]{Corollary}
\newtheorem{exAux}[theorem]{Example}

\newtheorem{Def}[theorem]{Definition}
\newenvironment{definition}{\begin{Def} \rm}{\end{Def}}
\newtheorem{Note}[theorem]{Note}
\newenvironment{note}{\begin{Note} \rm}{\end{Note}}
\newtheorem{Problem}[theorem]{Problem}

\newtheorem{Rem}[theorem]{Remark}

\newtheorem{Not}[theorem]{Notation}
\newenvironment{notation}{\begin{Not} \rm}{\end{Not}}
\newtheorem{Ass}[theorem]{Assumption}

\newenvironment{proof}{\medskip\noindent{\bf Proof.\ }}{\qed\medskip}

\newcommand{\qed}{\hfill\mbox{$\Box$\qquad\qquad}}

\newcommand{\Mat}[1]{\text{\rm Mat}_{#1}(\mathbb{K})}

\renewcommand{\b}[1]{\langle #1 \rangle}
\newcommand{\tr}{\text{\rm tr}}

\newcommand{\vphi}{\varphi}
\renewcommand{\th}{\theta}
\renewcommand{\indent}{\hspace{6mm}}

\begin{document}
\thispagestyle{empty}

\begin{center}
\LARGE \bf
\noindent
Transition maps between the 24 bases \\
for a Leonard pair
\end{center}

\smallskip

\begin{center}
\Large
Kazumasa Nomura and Paul Terwilliger
\end{center}

\smallskip

\begin{quote}
\small 
\begin{center}
\bf Abstract
\end{center}
\indent
Let $V$ denote a vector space with finite positive dimension.
We consider a pair of linear transformations 
$A : V \to V$ and $A^* : V \to V$ that satisfy (i) and (ii) below:
\begin{itemize}
\item[(i)] There exists a basis for $V$ with respect to which the
matrix representing $A$ is irreducible tridiagonal and the matrix
representing $A^*$ is diagonal.
\item[(ii)] There exists a basis for $V$ with respect to which the
matrix representing $A^*$ is irreducible tridiagonal and the matrix
representing $A$ is diagonal.
\end{itemize}
We call such a pair a {\em Leonard pair} on $V$.
In an earlier paper we described $24$ special bases for $V$. 
One feature of these bases is that with respect to each of them 
the matrices that represent $A$ and $A^*$ are 
(i) diagonal and irreducible tridiagonal
or (ii)  irreducible tridiagonal and diagonal
or (iii) lower bidiagonal and upper bidiagonal
or (iv) upper bidiagonal and lower bidiagonal.
For each ordered pair of bases among the $24$, there exists
a unique linear transformation from $V$ to $V$ that sends
the first basis to the second basis; we call this
the transition map. 
In this paper we find each transition map explicitly as a polynomial in $A,A^*$.
\end{quote}

\section{Leonard pairs}

\indent
We begin by recalling the notion of a Leonard pair.
We will use the following terms.
A square matrix $X$ is said to be {\em tridiagonal}
whenever each nonzero entry lies on either the diagonal, the subdiagonal,
or the superdiagonal. Assume $X$ is tridiagonal.
Then $X$ is said to be {\em irreducible}
whenever each entry on the subdiagonal is nonzero and each entry on
the superdiagonal is nonzero.
We now define a Leonard pair.
For the rest of this paper $\mathbb{K}$ will denote a field.

\begin{definition}  \cite{T:Leonard}  \label{def:LP}      \samepage
Let $V$ denote a vector space over $\mathbb{K}$ with finite positive
dimension.
By a {\em Leonard pair} on $V$ we mean an ordered pair $A,A^*$
where $A:V \to V$ and $A^*:V \to V$ are linear transformations
that satisfy (i) and (ii) below:
\begin{itemize}
\item[(i)] There exists a basis for $V$ with respect to which the
matrix representing $A$ is irreducible tridiagonal and the matrix
representing $A^*$ is diagonal.
\item[(ii)] There exists a basis for $V$ with respect to which the
matrix representing $A^*$ is irreducible tridiagonal and the matrix
representing $A$ is diagonal.
\end{itemize}
\end{definition}

\begin{note}
It is a common notational convention to use $A^*$ to represent the
conjugate-transpose of $A$. We are not using this convention.
In a Leonard pair $A,A^*$ the linear transformations $A$ and
$A^*$ are arbitrary subject to (i) and (ii) above.
\end{note}

We refer the reader to 
\cite{C:spinLP,H,N:aw,NT:balanced,NT:formula,NT:det,NT:mu,NT:span,NT:switch,
NT:affine,P,T:sub1,T:sub3,T:Leonard,T:24points,T:canform,T:intro,
T:intro2,T:split,T:array,T:qRacah,T:survey,TV,V,V:AW,ZK}
for background on Leonard pairs.
We especially recommend the survey \cite{T:survey}.
See \cite{AC,AC2,AC3,Bas,BT:Borel,BT:loop,Bow,C:mlt,F:RL,H:tetra,HT:tetra,
ITT,IT:shape,IT:uqsl2hat,IT:non-nilpotent,IT:tetra,IT:inverting,
IT:drg,IT:loop,IT:Krawt,ITW:equitable,M:LT,N:refine,N:height1,
NT:tde,R:multi,R:6j,T:qSerre,T:Kac-Moody,Vidar,Z}
for related topics.

\section{Leonard systems}

\indent
When working with a Leonard pair, it is convenient to consider a closely
related object called a {\em Leonard system}. 
To prepare for our definition
of a Leonard system, we recall a few concepts from linear algebra.
Let $d$ denote a nonnegative integer and let
$\Mat{d+1}$ denote the $\mathbb{K}$-algebra consisting of all $d+1$ by
$d+1$ matrices that have entries in $\mathbb{K}$. 
We index the rows and  columns by $0, 1, \ldots, d$. 
We let $\mathbb{K}^{d+1}$ denote the $\mathbb{K}$-vector space of all
$d+1$ by $1$ matrices that have entries in $\mathbb{K}$. We index the
rows by $0,1,\ldots,d$. We view $\mathbb{K}^{d+1}$ as a left module
for $\Mat{d+1}$. We observe this module is irreducible. 
For the rest of this paper, let $\cal A$ denote a $\mathbb{K}$-algebra 
isomorphic to $\Mat{d+1}$ and let $V$ denote an irreducible left $\cal A$-module.
We remark that $V$ is unique
up to isomorphism of $\cal A$-modules, and that $V$ has dimension $d+1$.
By a {\em basis} for $V$ we mean a sequence
of vectors that are linear independent and span $V$.
We emphasize that the ordering is important.
Let $\{v_i\}_{i=0}^d$ denote a basis for $V$.
For $X \in {\cal A}$ and $Y \in \Mat{d+1}$, we say 
{\em $Y$ represents $X$ with respect to} $\{v_i\}_{i=0}^d$
whenever $X v_j = \sum_{i=0}^d Y_{ij} v_i$ for $0 \leq j \leq d$.
For $A \in \cal A$ we say $A$ is {\em multiplicity-free}
whenever it has $d+1$ mutually distinct eigenvalues in $\mathbb{K}$. 
Assume $A$ is multiplicity-free. 
Let $\{\theta_i\}_{i=0}^d$ denote an ordering 
of the eigenvalues of $A$, and for $0 \leq i \leq d$ put
\begin{equation}        \label{eq:defEi}
    E_i = \prod_{\stackrel{0 \leq j \leq d}{j\neq i}}
             \frac{A-\theta_j I}{\theta_i - \theta_j},
\end{equation}
where $I$ denotes the identity of $\cal A$. 
We observe
(i) $AE_i = \theta_i E_i$ $(0 \leq i \leq d)$;
(ii) $E_i E_j = \delta_{i,j} E_i$ $(0 \leq i,j \leq d)$;
(iii) $\sum_{i=0}^{d} E_i = I$;
(iv) $A = \sum_{i=0}^{d} \theta_i E_i$.
Let $\cal D$ denote the subalgebra of $\cal A$ generated by $A$.
Using (i)--(iv) we find the sequence $\{E_i\}_{i=0}^d$
is a basis for the $\mathbb{K}$-vector space $\cal D$.
We call $E_i$ the {\em primitive idempotent} of $A$ associated with
$\theta_i$. 
It is helpful to think of these primitive idempotents as follows.
Observe 
\[
 V=E_0V+E_1V+\cdots+E_dV  \qquad \text{(direct sum)}.
\]
For $0 \leq i \leq d$, $E_iV$ is the (one dimensional) eigenspace of $A$
in $V$ associated with the eigenvalue $\theta_i$, and $E_i$ acts on $V$
as the projection onto this eigenspace.

\medskip

By a {\em Leonard pair in $\cal A$} we mean an ordered pair of elements
taken from $\cal A$ that act on $V$ as a Leonard pair in the sense of
Definition \ref{def:LP}.
We now define a Leonard system.

\medskip

\begin{definition}  \cite{T:Leonard}     \label{def:LS}   \samepage
By a {\em Leonard system} in $\cal A$ we mean a sequence
\[
  \Phi= (A; \{E_i\}_{i=0}^d; A^*; \{E^*_i\}_{i=0}^d)
\]
that satisfies (i)--(v) below.
\begin{itemize}
\item[(i)] Each of $A$, $A^*$ is a multiplicity-free element in $\cal A$.
\item[(ii)] $\{E_i\}_{i=0}^d$  is an ordering of the
   primitive idempotents of $A$.
\item[(iii)] $\{E^*_i\}_{i=0}^d$ is an ordering of the
   primitive idempotents of $A^*$.
\item[(iv)] For $0 \leq i,j \leq d$, 
\begin{equation}           \label{eq:Astrid}
   E_i A^* E_j =
    \begin{cases}  
        0 & \text{\rm if $|i-j|>1$},  \\
        \neq 0 & \text{\rm if $|i-j|=1$}.
    \end{cases}
\end{equation}
\item[(v)] For $0 \leq i,j \leq d$, 
\begin{equation}             \label{eq:Atrid}
   E^*_i A E^*_j =
    \begin{cases}  
        0 & \text{\rm if $|i-j|>1$},  \\
        \neq 0 & \text{\rm if $|i-j|=1$}.
    \end{cases}
\end{equation}
\end{itemize}
\end{definition}

\medskip

Leonard systems are related to Leonard pairs as follows.
Let $(A; \{E_i\}_{i=0}^d;A^*; \{E^*_i\}_{i=0}^d)$ denote a Leonard system
in $\cal A$. Then $A,A^*$ is a Leonard pair in $\cal A$
\cite[Section 3]{T:qRacah}.
Conversely, suppose $A,A^*$ is a Leonard pair in $\cal A$.
Then each of $A,A^*$ is multiplicity-free \cite[Lemma 1.3]{T:Leonard}.
Moreover there exists an ordering $\{E_i\}_{i=0}^d$ of the
primitive idempotents of $A$ and 
there exists an ordering $\{E^*_i\}_{i=0}^d$  of the
primitive idempotents of $A^*$ such that
$(A; \{E_i\}_{i=0}^d; A^*; \{E^*_i\}_{i=0}^d)$
is a Leonard system in $\cal A$ \cite[Lemma 3.3]{T:qRacah}.

\section{The 24 bases}

\indent
Let $\Phi=(A; \{E_i\}_{i=0}^d;A^*; \{E^*_i\}_{i=0}^d)$
denote a Leonard system in $\cal A$. In \cite{T:24points}
we obtained $24$ special bases for $V$, on which $A$ and $A^*$ act
in an attractive fashion. In this section we review these
bases. First we recall some notation.

\medskip

\begin{definition}        \label{def:th}       \samepage
Let $\Phi=(A; \{E_i\}_{i=0}^d;A^*; \{E^*_i\}_{i=0}^d)$ denote a 
Leonard system in $\cal A$.
For $0 \leq i \leq d$ we let $\theta_i$ (resp. $\theta^*_i$)
denote the eigenvalue of $A$ (resp. $A^*$) associated with
$E_i$ (resp. $E^*_i$).
We refer to $\{\theta_i\}_{i=0}^d$ (resp. $\{\theta^*_i\}_{i=0}^d$)
as the {\em eigenvalue sequence} (resp. {\em dual eigenvalue sequence})
of $\Phi$.
We observe $\{\theta_i\}_{i=0}^d$ (resp. $\{\theta^*_i\}_{i=0}^d$)
are mutually distinct
and contained in $\mathbb{K}$.
\end{definition}

\medskip

For an indeterminate $\lambda$ let $\mathbb{K}[\lambda]$ denote the 
$\mathbb{K}$-algebra of all polynomials in $\lambda$ that have
coefficients in $\mathbb{K}$.

\medskip

\begin{definition}             \label{def:tau}             \samepage
Referring to Definition \ref{def:th}, 
for $0 \leq i \leq d$ we define polynomials 
$\tau_i$, $\eta_i$, $\tau^*_i$, $\eta^*_i$ in $\mathbb{K}[\lambda]$
as follows:
\begin{align*}
 \tau_i &= (\lambda-\theta_0)(\lambda-\theta_1)\cdots(\lambda-\theta_{i-1}), \\
 \eta_i &= (\lambda-\theta_d)(\lambda-\theta_{d-1})\cdots(\lambda-\theta_{d-i+1}),\\
 \tau^*_i &= 
   (\lambda-\theta^*_0)(\lambda-\theta^*_1)\cdots(\lambda-\theta^*_{i-1}), \\
\eta^*_i &= 
   (\lambda-\theta^*_d)(\lambda-\theta^*_{d-1})\cdots(\lambda-\theta^*_{d-i+1}).
\end{align*}
Note that each of
$\tau_i$, $\eta_i$, $\tau^*_i$, $\eta^*_i$ is monic with degree $i$.
\end{definition}

\begin{lemma} {\rm \cite[Theorem 9.1]{T:24points}} \label{lem:24bases} \samepage
Let $\Phi=(A; \{E_i\}_{i=0}^d;A^*; \{E^*_i\}_{i=0}^d)$ denote a 
Leonard system in $\cal A$ and let
$\xi_0$, $\xi_d$, $\xi^*_0$, $\xi^*_d$ denote nonzero vectors in $V$ such that
\begin{equation}          \label{eq:defv0vd}
   \xi_0 \in E_0V,  \qquad
   \xi_d \in E_dV,  \qquad
   \xi^*_0 \in E^*_0V, \qquad
   \xi^*_d \in E^*_dV.
\end{equation}
Then each of the $24$ sequences below is a basis for $V$:
\[
\begin{array}{llll}
\{E_i \xi^*_0\}_{i=0}^d, &
\{E_i \xi^*_d\}_{i=0}^d, &
\{E_{d-i} \xi^*_0\}_{i=0}^d, &
\{E_{d-i} \xi^*_d\}_{i=0}^d, 
\\[7pt]
\{E^*_i \xi_0\}_{i=0}^d, &
\{E^*_i \xi_d\}_{i=0}^d, &
\{E^*_{d-i} \xi_0\}_{i=0}^d, &
\{E^*_{d-i} \xi_d\}_{i=0}^d, 
\\[7pt]
\{\tau_i(A)\xi^*_0\}_{i=0}^d, &
\{\tau_i(A)\xi^*_d\}_{i=0}^d, &
\{\eta_i(A)\xi^*_0\}_{i=0}^d, &
\{\eta_i(A)\xi^*_d\}_{i=0}^d, 
\\[7pt]
\{\tau^*_{d-i}(A^*)\xi_0\}_{i=0}^d, &
\{\tau^*_{d-i}(A^*)\xi_d\}_{i=0}^d, &
\{\eta^*_{d-i}(A^*)\xi_0\}_{i=0}^d, & 
\{\eta^*_{d-i}(A^*)\xi_d\}_{i=0}^d,
\\[7pt]
\{\tau_{d-i}(A)\xi^*_0\}_{i=0}^d, &
\{\tau_{d-i}(A)\xi^*_d\}_{i=0}^d, &
\{\eta_{d-i}(A)\xi^*_0\}_{i=0}^d, &  
\{\eta_{d-i}(A)\xi^*_d\}_{i=0}^d, 
\\[7pt]
\{\tau^*_i(A^*)\xi_0\}_{i=0}^d, &
\{\tau^*_i(A^*)\xi_d\}_{i=0}^d, &
\{\eta^*_i(A^*)\xi_0\}_{i=0}^d, & 
\{\eta^*_i(A^*)\xi_d\}_{i=0}^d.
\end{array}
\]
\end{lemma}

\begin{note}
Referring to Lemma 3.3 and with respect to each of the $24$ bases, 
the matrices that represent $A$ and $A^*$ are 
(i) diagonal and irreducible tridiagonal (row $1$);
(ii) irreducible tridiagonal and diagonal (row $2$);
(iii) lower bidiagonal and upper bidiagonal (rows $3,4$);
(iv) upper bidiagonal and lower bidiagonal (rows $5,6$).
See \cite{T:24points} for more information about the $24$ bases.
\end{note}

\medskip

Referring to Lemma \ref{lem:24bases}, for each ordered pair of bases
among the $24$ there exists a unique linear transformation 
from $V$ to $V$ that sends the first basis to the second basis; 
we call this the {\em transition map}.
In this paper we find each transition map explicitly as a polynomial in $A,A^*$.
Before we display the transition maps, we will review some basic
facts and prove a few lemmas about Leonard pairs.

\section{The $D_4$ action}

\indent
Let $\Phi=(A; \{E_i\}_{i=0}^d; A^*; \{E^*_i\}_{i=0}^d)$
denote a Leonard system in $\cal A$.
Then each of the following is a Leonard system in $\cal A$:
\begin{align*}
\Phi^{*}  &:= 
       (A^*; \{E^*_i\}_{i=0}^d; A; \{E_i\}_{i=0}^d), \\
\Phi^{\downarrow} &:=
       (A; \{E_i\}_{i=0}^d; A^*; \{E^*_{d-i}\}_{i=0}^d), \\
\Phi^{\Downarrow} &:=
       (A; \{E_{d-i}\}_{i=0}^d; A^*; \{E^*_{i}\}_{i=0}^d).
\end{align*}
Viewing $*$, $\downarrow$, $\Downarrow$ as permutations on the set of
all the Leonard systems,
\begin{equation}    \label{eq:relation1}
\text{$*^2$$\;=\;$$\downarrow^2$$\;=\;$$\Downarrow^2$$\;=\;$$1$},
\end{equation}
\begin{equation}    \label{eq:relation2}
\text{$\Downarrow$$*$$\;=\;$$*$$\downarrow$}, \quad
\text{$\downarrow$$*$$\;=\;$$*$$\Downarrow$}, \quad
\text{$\downarrow$$\Downarrow$$\;=\;$$\Downarrow$$\downarrow$}.
\end{equation}
The group generated by symbols $*$, $\downarrow$, $\Downarrow$ subject
to the relations (\ref{eq:relation1}), (\ref{eq:relation2}) is the
dihedral group $D_4$. We recall $D_4$ is the group of symmetries of a
square, and has $8$ elements.
Apparently $*$, $\downarrow$, $\Downarrow$ induce an action of $D_4$
on the set of all Leonard systems.
Two Leonard systems will be called {\em relatives} whenever they are
in the same orbit of this $D_4$ action. 
The relatives of $\Phi$ are as follows:

\medskip
\noindent
\begin{center}
\begin{tabular}{c|c}
name  &  relative \\
\hline
$\Phi$ & 
       $(A; \{E_i\}_{i=0}^d; A^*;  \{E^*_i\}_{i=0}^d)$ \\ 
$\Phi^{\downarrow}$ &
       $(A; \{E_i\}_{i=0}^d; A^*;  \{E^*_{d-i}\}_{i=0}^d)$ \\ 
$\Phi^{\Downarrow}$ &
       $(A; \{E_{d-i}\}_{i=0}^d; A^*;  \{E^*_i\}_{i=0}^d)$ \\ 
$\Phi^{\downarrow \Downarrow}$ &
       $(A; \{E_{d-i}\}_{i=0}^d; A^*;  \{E^*_{d-i}\}_{i=0}^d)$ \\ 
$\Phi^{*}$  & 
       $(A^*; \{E^*_i\}_{i=0}^d; A;  \{E_i\}_{i=0}^d)$ \\ 
$\Phi^{\downarrow *}$ &
       $(A^*; \{E^*_{d-i}\}_{i=0}^d; A;  \{E_i\}_{i=0}^d)$ \\ 
$\Phi^{\Downarrow *}$ &
       $(A^*; \{E^*_i\}_{i=0}^d; A;  \{E_{d-i}\}_{i=0}^d)$ \\ 
$\Phi^{\downarrow \Downarrow *}$ &
       $(A^*; \{E^*_{d-i}\}_{i=0}^d; A;  \{E_{d-i}\}_{i=0}^d)$
\end{tabular}
\end{center}

\section{The parameter array}

\indent
In this section we recall the parameter array of a Leonard system.

\medskip

\begin{definition}  \cite[Theorem 4.6]{NT:formula}  \label{def:splitseq} \samepage
Let $\Phi=(A; \{E_i\}_{i=0}^d;A^*; \{E^*_i\}_{i=0}^d)$ denote a 
Leonard system with eigenvalue sequence $\{\theta_i\}_{i=0}^d$
and dual eigenvalue sequence $\{\theta^*_i\}_{i=0}^d$.
For $1 \leq i \leq d$ we define scalars
\begin{align}
 \varphi_i 
   &= (\theta^*_0-\theta^*_i)
       \frac{\text{tr}(\tau_{i}(A)E^*_0)}
            {\text{tr}(\tau_{i-1}(A)E^*_0)},
                               \label{eq:defvarphi}   \\
 \phi_i
   &= (\theta^*_0-\theta^*_i)
       \frac{\text{tr}(\eta_i(A)E^*_0)}
            {\text{tr}(\eta_{i-1}(A)E^*_0)}
                               \label{eq:defphi},
\end{align}
where tr means trace.
In (\ref{eq:defvarphi}), (\ref{eq:defphi}) the denominators 
are nonzero by \cite[Corollary 4.5]{NT:formula}.
Also by \cite[Corollary 4.5]{NT:formula}
each of $\varphi_i$, $\phi_i$ is nonzero
for $1 \leq i \leq d$.
The sequence $\{\varphi_i\}_{i=1}^d$ 
(resp. $\{\phi_i\}_{i=1}^d$) is called the 
{\em first split sequence} (resp. {\em second split sequence}) of $\Phi$.
\end{definition}

\begin{definition}                              \samepage
Let $\Phi=(A; \{E_i\}_{i=0}^d;A^*; \{E^*_i\}_{i=0}^d)$ denote a 
Leonard system.
By the {\em parameter array} of $\Phi$ we mean the sequence
\begin{equation}      \label{eq:param}
      (\{\theta_i\}_{i=0}^d; \{\theta^*_i\}_{i=0}^d;
        \{\varphi_i\}_{i=1}^d; \{\phi_i\}_{i=1}^d),
\end{equation}
where $\{\theta_i\}_{i=0}^d$ (resp. $\{\theta^*_i\}_{i=0}^d$)
denotes the eigenvalue sequence (resp. dual eigenvalue sequence) of $\Phi$
and $\{\varphi_i\}_{i=1}^d$ (resp. $\{\phi_i\}_{i=1}^d$)
denotes the first split sequence (resp. second split sequence) of $\Phi$.
For notational convenience we abbreviate
\[
  \varphi:=\varphi_1\varphi_2\cdots\varphi_d,  \qquad\qquad
  \phi:=\phi_1\phi_2\cdots\phi_d.
\]
\end{definition}

\medskip

The $D_4$ action affects the parameter array as follows.

\medskip

\begin{lemma} {\rm \cite[Theorem 1.11]{T:Leonard}}  \label{lem:D4} \samepage
Let $\Phi=(A; \{E_i\}_{i=0}^d; A^*; \{E^*_i\}_{i=0}^d)$ denote
a Leonard system  and let 
$(\{\theta_i\}_{i=0}^d; \{\theta^*_i\}_{i=0}^d;
  \{\varphi_i\}_{i=1}^d; \{\phi_i\}_{i=1}^d)$
denote the parameter array of $\Phi$.
Then the following (i)--(iii) hold.
\begin{itemize}
\item[(i)] The parameter array of $\Phi^*$ is 
\[
   (\{\theta^*_i\}_{i=0}^d; \{\theta_i\}_{i=0}^d;
        \{\varphi_i\}_{i=1}^d; \{\phi_{d-i+1}\}_{i=1}^d).
\]
\item[(ii)]
The parameter array of $\Phi^{\downarrow}$ is
\[
  (\{\theta_i\}_{i=0}^d; \{\theta^*_{d-i}\}_{i=0}^d;
        \{\phi_{d-i+1}\}_{i=1}^d; \{\varphi_{d-i+1}\}_{i=1}^d).
\]
\item[(iii)]
The parameter array of $\Phi^{\Downarrow}$ is
\[
  (\{\theta_{d-i}\}_{i=0}^d; \{\theta^*_{i}\}_{i=0}^d;
        \{\phi_{i}\}_{i=1}^d; \{\varphi_{i}\}_{i=1}^d).
\]
\end{itemize}
\end{lemma}

\section{The antiautomorphism $\dagger$}

\indent
Associated with a given Leonard system in $\cal A$, there is a certain
antiautomorphism of $\cal A$ denoted by $\dagger$ and defined below.
Recall an {\em antiautomorphism} of $\cal A$ is an isomorphism
of $\mathbb{K}$-vector spaces $\sigma : {\cal A} \to {\cal A}$ 
such that  $(XY)^\sigma = Y^\sigma X^\sigma$
for all $X,Y \in {\cal A}$.

\medskip

\begin{lemma} {\rm \cite[Theorem 6.1]{T:qRacah}} \label{lem:dagger} \samepage
Let $(A; \{E_i\}_{i=0}^d; A^*; \{E^*_i\}_{i=0}^d)$
denote a Leonard system in $\cal A$.
Then there exists a unique antiautomorphism $\dagger$ of $\cal A$ such that
$A^\dagger=A$ and $A^{*\dagger}=A^*$.
\end{lemma}

\begin{lemma} {\rm \cite[Lemma 6.3]{T:qRacah}}  \label{lem:Eidagger} \samepage
Let $(A; \{E_i\}_{i=0}^d; A^*; \{E^*_i\}_{i=0}^d)$
denote a Leonard system in $\cal A$ and let $\dagger$ denote
the corresponding antiautomorphism from Lemma \ref{lem:dagger}.
Then (i), (ii) hold below.
\begin{itemize}
\item[(i)]
Let $\cal D$ denote the subalgebra of $\cal A$ generated by $A$.
Then $X^\dagger=X$ for all $X \in {\cal D}$;
in particular $E^\dagger_i=E_i$ for $0 \leq i \leq d$.
\item[(ii)]
Let ${\cal D}^*$ denote the subalgebra of $\cal A$ generated by $A^*$.
Then $X^\dagger=X$ for all $X \in {\cal D}^*$;
in particular $E^{*\dagger}_i=E^*_i$ for $0 \leq i \leq d$.
\end{itemize}
\end{lemma}

\section{Some reduction rules}

\indent
In this section we give some formulae involving Leonard systems which 
we will use later in the paper.

\medskip

\begin{lemma}                  \label{lem:tauiA}       \samepage
Let $\Phi=(A; \{E_i\}_{i=0}^d; A^*; \{E^*_i\}_{i=0}^d)$ denote
a Leonard system and let 
$(\{\theta_i\}_{i=0}^d$; 
$\{\theta^*_i\}_{i=0}^d; \{\varphi_i\}_{i=1}^d; \{\phi_i\}_{i=1}^d)$
denote the parameter array of $\Phi$.
Then for $0 \leq i,j \leq d$ we have
\begin{align}    
 E_0 \tau^*_i(A^*)\tau_j(A)E^*_0 &= 
  \delta_{i,j}\varphi_1\varphi_2\cdots\varphi_i E_0E^*_0, 
                                                      \label{eq:basic} \\[5pt]
 E_0 \eta^*_i(A^*)\tau_j(A)E^*_d &= 
  \delta_{i,j}\phi_d\phi_{d-1}\cdots\phi_{d-i+1} E_0E^*_d, 
                                                      \label{eq:basicd} \\[5pt]
 E_d \tau^*_i(A^*)\eta_j(A)E^*_0 &= 
  \delta_{i,j}\phi_1\phi_2\cdots\phi_i E_dE^*_0, 
                                                      \label{eq:basicD} \\[5pt]
 E_d \eta^*_i(A^*)\eta_j(A)E^*_d &= 
  \delta_{i,j}\vphi_d\vphi_{d-1}\cdots\vphi_{d-i+1} E_dE^*_d, 
                                                      \label{eq:basicdD} \\[5pt]
 E^*_0 \tau_i(A)\tau^*_j(A^*)E_0 &= 
  \delta_{i,j}\varphi_1\varphi_2\cdots\varphi_i E^*_0E_0, 
                                                      \label{eq:basics} \\[5pt]
 E^*_0 \eta_i(A)\tau^*_j(A^*)E_d &= 
  \delta_{i,j}\phi_1\phi_2\cdots\phi_i E^*_0E_d, 
                                                      \label{eq:basicsds} \\[5pt]
 E^*_d \tau_i(A)\eta^*_j(A^*)E_0 &= 
  \delta_{i,j}\phi_d\phi_{d-1}\cdots\phi_{d-i+1} E^*_dE_0, 
                                                      \label{eq:basicsDs} \\[5pt]
 E^*_d \eta_i(A)\eta^*_j(A^*)E_d &= 
  \delta_{i,j}\vphi_d\vphi_{d-1}\cdots\vphi_{d-i+1} E^*_dE_d.
                                                      \label{eq:basicdDs} 
\end{align}
\end{lemma}

\begin{proof}
We first show \eqref{eq:basic}.
Replacing $\Phi$ by $\Phi^*$ and using $\dagger$ if necessary 
we may assume $i \geq j$.
In the left-hand side of \eqref{eq:basic} we insert a factor $I$ between
$\tau^*_i(A^*)$ and $\tau_j(A)$. 
We expand using $I=\sum_{r=0}^d E^*_r$ and simplify
the result using $\tau^*_i(A^*)E^*_r= \tau^*_i(\theta^*_r)E^*_r$ for
$0 \leq r \leq d$.
By these comments the left-hand side of \eqref{eq:basic} is equal to
\begin{equation}              \label{eq:basicaux}
 \sum_{r=0}^d \tau^*_i(\theta^*_r) E_0 E^*_r \tau_j(A)E^*_0.
\end{equation}
For $0 \leq r \leq d$ we examine term $r$ in \eqref{eq:basicaux}.
By Definition \ref{def:tau} we have $\tau^*_i(\theta^*_r)= 0$  if $r<i$.
By \cite[Lemma 5.10(i)]{T:qRacah} $E^*_rA^sE^*_0= 0$ for $0 \leq s \leq r-1$.
By this and since the polynomial $\tau_j$ is  monic of degree $j$,
$E^*_r \tau_j(A)E^*_0$ is $0$ if $j<r$ and  $E^*_r A^r E^*_0$  if $j=r$.
Referring to \cite[Theorem 8.4]{T:qRacah} we have $E^*_iA^iE^*_0= p_i(A)E^*_0$.
Note that $E_0p_i(A)= p_i(\theta_0)E_0$ by construction and
$\tau^*_i(\theta^*_i)p_i(\theta_0)= \varphi_1 \varphi_2 \cdots \varphi_i$
by \cite[Lemma 17.5]{T:qRacah}.
By these comments and since $i\geq j$
the sum \eqref{eq:basicaux} is equal to the right-hand side of \eqref{eq:basic}.
We have now verified \eqref{eq:basic}.
To get \eqref{eq:basicd}--\eqref{eq:basicdDs}, apply $D_4$ to \eqref{eq:basic}
and use Lemma \ref{lem:D4}.
\end{proof}

\begin{lemma}        \label{lem:simplify} \samepage
Let $\Phi=(A; \{E_i\}_{i=0}^d; A^*; \{E^*_i\}_{i=0}^d)$ denote
a Leonard system and let 
$(\{\theta_i\}_{i=0}^d$; 
$\{\theta^*_i\}_{i=0}^d;   \{\varphi_i\}_{i=1}^d; \{\phi_i\}_{i=1}^d)$
denote the parameter array of $\Phi$.
Then for $0 \leq r \leq d$ we have
\begin{align}  
 E^*_0E_dE^*_dE_r &=
   \frac{\vphi}{\tau_d(\th_d)\tau^*_d(\th^*_d)}
   \frac{\phi_d\phi_{d-1}\cdots\phi_{d-r+1}}
        {\vphi_1\vphi_2\cdots\vphi_r}
    E^*_0E_r,                             \label{eq:Es0EdEsdEr}  \\
 E^*_0E_0E^*_dE_r &=
   \frac{\vphi}{\eta_d(\th_0)\tau^*_d(\th^*_d)}
   \frac{\phi_d\phi_{d-1}\cdots\phi_{d-r+1}}
        {\vphi_1\vphi_2\cdots\vphi_r}
    E^*_0E_r,                             \label{eq:Es0E0EsdEr}  \\
 E^*_dE_dE^*_0E_r &=
   \frac{\phi}{\tau_d(\th_d)\eta^*_d(\th^*_0)}
   \frac{\vphi_1\vphi_2\cdots\vphi_r}
        {\phi_d\phi_{d-1}\cdots\phi_{d-r+1}}  
   E^*_dE_r,                             \label{eq:EsdEdEs0Er}    \\
 E^*_dE_0E^*_0E_r &=
   \frac{\phi}{\eta_d(\th_0)\eta^*_d(\th^*_0)}
   \frac{\vphi_1\vphi_2\cdots\vphi_r}
        {\phi_d\phi_{d-1}\cdots\phi_{d-r+1}}
    E^*_dE_r,                             \label{eq:EsdE0Es0Er}  \\
 E_0E^*_dE_dE^*_r &=
   \frac{\vphi}{\tau_d(\th_d)\tau^*_d(\th^*_d)}
   \frac{\phi_1\phi_2\cdots\phi_r}
        {\vphi_1\vphi_2\cdots\vphi_r}
    E_0E^*_r,                             \label{eq:E0EsdEdEsr}  \\
 E_0E^*_0E_dE^*_r &=
   \frac{\vphi}{\tau_d(\th_d)\eta^*_d(\th^*_0)}
   \frac{\phi_1\phi_2\cdots\phi_r}
        {\vphi_1\vphi_2\cdots\vphi_r}
    E_0E^*_r,                             \label{eq:E0Es0EdEsr}  \\
 E_dE^*_dE_0E^*_r &=
   \frac{\phi}{\eta_d(\th_0)\tau^*_d(\th^*_d)}
   \frac{\vphi_1\vphi_2\cdots\vphi_r}
        {\phi_1\phi_2\cdots\phi_r}  
   E_dE^*_r,                             \label{eq:EdEsdE0Esr}  \\
 E_dE^*_0E_0E^*_r &=
   \frac{\phi}{\eta_d(\th_0)\eta^*_d(\th^*_0)}
   \frac{\vphi_1\vphi_2\cdots\vphi_r}
        {\phi_1\phi_2\cdots\phi_r}
    E_dE^*_r.                             \label{eq:EdEs0E0Esr} 
\end{align}
\end{lemma}

\begin{proof}
We first show \eqref{eq:Es0EdEsdEr}.
Following \cite[Definition 5.1]{NT:switch} we define
\begin{equation}                \label{eq:defS}
 S = \sum_{i=0}^d 
       \frac{\phi_d\phi_{d-1}\cdots\phi_{d-i+1}}
            {\vphi_1\vphi_2\cdots\vphi_i}   E_i.
\end{equation}
By \cite[Lemma 11.5]{NT:switch},
\begin{equation}          \label{eq:simplifyaux1}
 SE^*_0 =
  \frac{\tau_d(\th_d)\tau^*_d(\th^*_d)}
       {\vphi}
  E^*_dE_dE^*_0.
\end{equation}
Applying $\dagger$ to \eqref{eq:simplifyaux1} using Lemma \ref{lem:Eidagger},
and multiplying the result on the right by $E_r$ we find
\begin{equation}          \label{eq:simplifyaux3}
 E^*_0SE_r =
  \frac{\tau_d(\theta_d)\tau^*_d(\theta^*_d)}
       {\varphi}
  E^*_0E_dE^*_dE_r.
\end{equation}
Using \eqref{eq:defS} we find
\begin{equation}            \label{eq:simplifyaux4}
  SE_r =
    \frac{\phi_d\phi_{d-1}\cdots\phi_{d-r+1}}
         {\varphi_1\varphi_2\cdots\varphi_r} E_r.
\end{equation}
Combining \eqref{eq:simplifyaux3} and \eqref{eq:simplifyaux4} we obtain
\eqref{eq:Es0EdEsdEr}.
To get \eqref{eq:Es0E0EsdEr}--\eqref{eq:EdEs0E0Esr},
apply $D_4$ to \eqref{eq:Es0EdEsdEr} and use Lemma \ref{lem:D4}.
\end{proof}

\section{Some traces}

\indent
Let $\Phi=(A; \{E_i\}_{i=0}^d;A^*; \{E^*_i\}_{i=0}^d)$ denote a 
Leonard system.
In our main results we will need the scalars
\[
 \tr(E_rE^*_0),\qquad \tr(E_rE^*_d), \qquad \tr(E^*_rE_0), \qquad \tr(E^*_rE_d)
\]
for $0 \leq r \leq d$.
In this section we make some comments about these scalars.

\medskip

\begin{lemma}  {\rm \cite[Lemma 9.2]{T:qRacah}}  \label{lem:mi} \samepage
Let $\Phi=(A; \{E_i\}_{i=0}^d;A^*; \{E^*_i\}_{i=0}^d)$ denote a 
Leonard system.
Then
\begin{align}
E_rE^*_0E_r &= \tr(E_rE^*_0)E_r, & 
E_rE^*_dE_r &= \tr(E_rE^*_d)E_r,       \label{eq:ErEs0Er} \\
E^*_rE_0E^*_r &= \tr(E^*_rE_0)E^*_r, & 
E^*_rE_dE^*_r &= \tr(E^*_rE_d)E^*_r,   \label{eq:EsrE0Esr} \\
E^*_0E_rE^*_0 &= \tr(E_rE^*_0)E^*_0, & 
E^*_dE_rE^*_d &= \tr(E_rE^*_d)E^*_d,   \label{eq:EsdErEsd} \\
E_0E^*_rE_0 &= \tr(E^*_rE_0)E_0, & 
E_dE^*_rE_d &= \tr(E^*_rE_d)E_d.     \label{eq:E0EsrE0}
\end{align}
\end{lemma}

\begin{lemma} {\rm \cite[Theorem 17.12]{T:qRacah}} \label{lem:trErEs0} \samepage
Let $\Phi=(A; \{E_i\}_{i=0}^d;A^*; \{E^*_i\}_{i=0}^d)$ denote a 
Leonard system and let
$ (\{\theta_i\}_{i=0}^d$; $\{\theta^*_i\}_{i=0}^d;
        \{\varphi_i\}_{i=1}^d; \{\phi_i\}_{i=1}^d)$
denote the parameter array of $\Phi$.
Then for $0 \leq r \leq d$ we have
\begin{align}     
 \tr(E_rE^*_0) &=
  \frac{\varphi_1\varphi_2\cdots\varphi_r\,\phi_1\phi_2\cdots\phi_{d-r}}
       {\eta^*_d(\theta^*_0)\tau_r(\theta_r)\eta_{d-r}(\theta_r)},
                      \label{eq:trErEs0}  \\
 \tr(E_rE^*_d) &=
  \frac{\phi_d\phi_{d-1}\cdots\phi_{d-r+1}\,
             \varphi_d\varphi_{d-1}\cdots\varphi_{r+1}}
       {\tau^*_d(\theta^*_d)\tau_r(\theta_r)\eta_{d-r}(\theta_r)},
                      \label{eq:trErEsd}  \\
 \tr(E^*_rE_0) &=
  \frac{\varphi_1\varphi_2\cdots\varphi_r\,
            \phi_d\varphi_{d-1}\cdots\phi_{r+1}}
       {\eta_d(\theta_0)\tau^*_r(\theta^*_r)\eta^*_{d-r}(\theta^*_r)},
                      \label{eq:trEsrE0}  \\
 \tr(E^*_rE_d) &=
  \frac{\phi_1\phi_2\cdots\phi_r\,
            \varphi_d\varphi_{d-1}\cdots\varphi_{r+1}}
       {\tau_d(\theta_d)\tau^*_r(\theta^*_r)\eta^*_{d-r}(\theta^*_r)}.
                      \label{eq:trEsrEd} 
\end{align}
\end{lemma}

\begin{corollary}       \label{cor:trE0Es0}        \samepage
Let $\Phi=(A; \{E_i\}_{i=0}^d;A^*; \{E^*_i\}_{i=0}^d)$ denote a 
Leonard system and let
$ (\{\theta_i\}_{i=0}^d$; $\{\theta^*_i\}_{i=0}^d;
        \{\varphi_i\}_{i=1}^d; \{\phi_i\}_{i=1}^d)$
denote the parameter array of $\Phi$.
Then for $0 \leq r \leq d$ we have
\begin{align}
 \tr(E_0E^*_0) &= 
  \frac{\phi}{\eta_d(\theta_0)\eta^*_d(\theta^*_0)},  &
 \tr(E_0E^*_d) &=
  \frac{\varphi}{\eta_d(\theta_0)\tau^*_d(\theta^*_d)},  \label{eq:trE0Es0}  \\
 \tr(E_dE^*_0) &= 
  \frac{\varphi}{\tau_d(\theta_d)\eta^*_d(\theta^*_0)},  &
 \tr(E_dE^*_d) &= 
  \frac{\phi}{\tau_d(\theta_d)\tau^*_d(\theta^*_d)}.   \label{eq:trEdEs0}
\end{align}
\end{corollary}

\begin{proof}
Follows from Lemma \ref{lem:trErEs0}.
\end{proof}

\begin{corollary}
Let $\Phi=(A; \{E_i\}_{i=0}^d;A^*; \{E^*_i\}_{i=0}^d)$ denote a 
Leonard system. 
Then each of
$\tr(E_rE^*_0)$,
$\tr(E_rE^*_d)$,
$\tr(E^*_rE_0)$,
$\tr(E^*_rE_d)$
is nonzero for $0 \leq r \leq d$.
\end{corollary}

\begin{proof}
By Lemma \ref{lem:trErEs0} and since each of $\vphi_i, \phi_i$ is
nonzero for $1 \leq i \leq d$.
\end{proof}

\section{A bilinear form}

\indent
In this section we associate with each Leonard system a certain
bilinear form.
To prepare for this we recall a few concepts from linear algebra.

By a {\em bilinear form on $V$} we mean a map
$\b{\;,\,} : V \times V \to \mathbb{K}$ that satisfies the following
four conditions for all $u,v,w \in V$ and for all $\alpha \in \mathbb{K}$:
(i) $\b{u+v,w}=\b{u,w}+\b{v,w}$;
(ii) $\b{\alpha u,v}=\alpha \b{u,v}$;
(iii) $\b{u,v+w}=\b{u,v}+\b{u,w}$;
(iv) $\b{u,\alpha v}=\alpha \b{u,v}$.
Let $\b{\;,\,}$ denote a bilinear form on $V$.
This form is said to be {\em symmetric} whenever 
$\b{u,v}=\b{v,u}$ for all $u,v \in V$.
Let $\b{\;,\,}$ denote a bilinear form on $V$.
Then the following are equivalent:
(i) there exists a nonzero $u \in V$ such that $\b{u,v}=0$ for all $v\in V$;
(ii) there exists a nonzero $v\in V$ such that $\b{u,v}=0$ for all $u\in V$.
The form $\b{\;,\,}$ is said to be {\em degenerate} whenever (i), (ii) hold
and {\em nondegenerate} otherwise.
Let $\sigma : {\cal A} \to {\cal A}$ denote an antiautomorphism.
Then there exists a nonzero bilinear form $\b{\;,\,}$ on $V$
such that $\b{Xu,v}=\b{u,X^\sigma v}$ for all $u,v \in V$ and 
for all $X \in {\cal A}$.
The form is unique up to multiplication by a nonzero scalar in $\mathbb{K}$.
The form is nondegenerate.
We refer to this form as the {\em bilinear form on $V$ associated with $\sigma$}.
This form is not symmetric in general.

We now return our attention to Leonard systems.

\medskip

\begin{definition}           \label{def:bilin}        \samepage
Let $\Phi=(A; \{E_i\}_{i=0}^d; A^*; \{E^*_i\}_{i=0}^d)$
denote a Leonard system in $\cal A$.
Let $\dagger : {\cal A} \to {\cal A}$ denote the corresponding
antiautomorphism from Lemma \ref{lem:dagger}.
For the rest of this paper we let $\b{\;,\,}$ denote the bilinear form
on $V$ associated with $\dagger$.
By construction, for $X \in {\cal A}$ we have
\begin{equation}            \label{eq:genbilin}
  \b{Xu,v}=\b{u,X^\dagger v}  \qquad\qquad (u,v \in V).
\end{equation}
\end{definition}

\medskip

\begin{lemma}  {\rm \cite[Lemma 15.2]{T:qRacah}}  \label{lem:bilin} \samepage
Let $\Phi=(A; \{E_i\}_{i=0}^d; A^*; \{E^*_i\}_{i=0}^d)$
denote a Leonard system in $\cal A$.
Let ${\cal D}$ (resp. ${\cal D}^*$) denote the subalgebra of $\cal A$
generated by $A$ (resp. $A^*$).
Then for $X \in {\cal D} \cup {\cal D}^*$ we have
\begin{equation}           \label{eq:bilin}
  \b{Xu,v} = \b{u,Xv}  \qquad\qquad (u,v \in V).
\end{equation}
\end{lemma}

\begin{lemma}  {\rm \cite[Corollary 15.4]{T:qRacah}}     \samepage
Let $\Phi=(A; \{E_i\}_{i=0}^d; A^*; \{E^*_i\}_{i=0}^d)$
denote a Leonard system in $\cal A$. 
Then the bilinear form $\b{\;,\,}$ is symmetric.
\end{lemma}

\medskip
Let $\Phi=(A; \{E_i\}_{i=0}^d; A^*; \{E^*_i\}_{i=0}^d)$
denote a Leonard system in $\cal A$.
In our discussion of $\Phi$ we will refer to the inner products
\begin{equation}             \label{eq:8products}
\begin{array}{ccccccc}
\b{\xi_0,\xi^*_0}, & \; &
\b{\xi_0,\xi^*_d}, & \; &
\b{\xi_d,\xi^*_0}, & \; &
\b{\xi_d,\xi^*_d}, \\
\b{\xi_0,\xi_0},  & \; &
\b{\xi_d,\xi_d}, &  \; &
\b{\xi^*_0,\xi^*_0}, & \; &
\b{\xi^*_d,\xi^*_d},
\end{array}
\end{equation}
where the vectors
$\xi_0,\xi_d,\xi^*_0,\xi^*_d$ are from Lemma \ref{lem:24bases}.
We have some comments on these inner products.

\medskip

\begin{lemma} {\rm \cite[Lemma 15.5]{T:qRacah}} \label{lem:nonzero}  \samepage
Let $\Phi=(A; \{E_i\}_{i=0}^d; A^*; \{E^*_i\}_{i=0}^d)$
denote a Leonard system in $\cal A$ and let
$\xi_0,\xi_d,\xi^*_0,\xi^*_d$ denote nonzero vectors in $V$
that satisfy \eqref{eq:defv0vd}.
Then each of the scalars \eqref{eq:8products} is nonzero.
\end{lemma}

\begin{lemma}  {\rm \cite[Lemma 15.5]{T:qRacah}}         \samepage
Let $\Phi=(A; \{E_i\}_{i=0}^d; A^*; \{E^*_i\}_{i=0}^d)$
denote a Leonard system in $\cal A$ and let
$\xi_0,\xi_d,\xi^*_0,\xi^*_d$ denote nonzero vectors in $V$
that satisfy \eqref{eq:defv0vd}. Then
\begin{align}
 E_0\xi^*_0 &= \frac{\b{\xi_0,\xi^*_0}}{\b{\xi_0,\xi_0}} \xi_0,   &
 E_d\xi^*_0 &= \frac{\b{\xi_d,\xi^*_0}}{\b{\xi_d,\xi_d}} \xi_d,
                                                           \label{eq:vs0} \\
 E_0\xi^*_d &= \frac{\b{\xi_0,\xi^*_d}}{\b{\xi_0,\xi_0}} \xi_0,   &
 E_d\xi^*_d &= \frac{\b{\xi_d,\xi^*_d}}{\b{\xi_d,\xi_d}} \xi_d,    
                                                           \label{eq:vsd} \\
 E^*_0\xi_0 &= \frac{\b{\xi_0,\xi^*_0}}{\b{\xi^*_0,\xi^*_0}} \xi^*_0,  &
 E^*_d\xi_0 &= \frac{\b{\xi_0,\xi^*_d}}{\b{\xi^*_d,\xi^*_d}} \xi^*_d,  
                                                             \label{eq:v0} \\
 E^*_0\xi_d &= \frac{\b{\xi_d,\xi^*_0}}{\b{\xi^*_0,\xi^*_0}} \xi^*_0, &
 E^*_d\xi_d &= \frac{\b{\xi_d,\xi^*_d}}{\b{\xi^*_d,\xi^*_d}} \xi^*_d.  
                                                              \label{eq:vd}
\end{align}
\end{lemma}

\medskip

In the following two lemmas we give some relationships involving 
the eight scalars in \eqref{eq:8products}.

\medskip

\begin{lemma}  {\rm \cite[Lemma 15.5]{T:qRacah}} \label{lem:bv0vs0} \samepage
Let $\Phi=(A; \{E_i\}_{i=0}^d; A^*; \{E^*_i\}_{i=0}^d)$
denote a Leonard system in $\cal A$ and let
$\xi_0,\xi_d,\xi^*_0,\xi^*_d$ denote nonzero vectors in $V$
that satisfy \eqref{eq:defv0vd}.
Then
\begin{align}
\b{\xi_0,\xi^*_0}^2 
 &= \text{\rm tr}(E_0E^*_0) \b{\xi_0,\xi_0}\b{\xi^*_0,\xi^*_0}, \label{eq:00s}  \\
\b{\xi_0,\xi^*_d}^2 
 &= \text{\rm tr}(E_0E^*_d) \b{\xi_0,\xi_0}\b{\xi^*_d,\xi^*_d}, \label{eq:0ds}  \\
\b{\xi_d,\xi^*_0}^2 
 &= \text{\rm tr}(E_dE^*_0) \b{\xi_d,\xi_d}\b{\xi^*_0,\xi^*_0}, \label{eq:d0s}  \\
\b{\xi_d,\xi^*_d}^2 
 &= \text{\rm tr}(E_dE^*_d) \b{\xi_d,\xi_d}\b{\xi^*_d,\xi^*_d}. \label{eq:dds}
\end{align}
\end{lemma}

\begin{lemma}      \label{lem:newrel}            \samepage
Let $\Phi=(A; \{E_i\}_{i=0}^d; A^*; \{E^*_i\}_{i=0}^d)$
denote a Leonard system in $\cal A$ and let
$\xi_0,\xi_d,\xi^*_0,\xi^*_d$ denote nonzero vectors in $V$
that satisfy \eqref{eq:defv0vd}.
Then
\begin{equation}         \label{eq:newrel}
  \frac{\b{\xi_0,\xi^*_0}}{\b{\xi_0,\xi^*_d}}
= \frac{\phi}{\varphi} \frac{\b{\xi_d,\xi^*_0}}{\b{\xi_d,\xi^*_d}}.
\end{equation}
\end{lemma}

\begin{proof}
Let $S$ be from \eqref{eq:defS}.
Recall $\xi_0 \in E_0V$ so $S \xi_0=\xi_0$.
Similarly $\xi_d \in E_dV$ so $S \xi_d= \phi\varphi^{-1} \xi_d$.
By \cite[Theorem 6.7]{NT:switch} we have $S E^*_0V=E^*_dV$
so there exists a nonzero $t \in \mathbb{K}$ such that
$S \xi^*_0 = t \xi^*_d$.
We may now argue
\[
 \b{\xi_d, \xi^*_d}
   = t^{-1}\b{\xi_d, S \xi^*_0}
   = t^{-1}\b{S\xi_d, \xi^*_0} 
   = \phi\varphi^{-1} t^{-1} \b{\xi_d, \xi^*_0}
\]
and
\[
 \b{\xi_0, \xi^*_d}
  = t^{-1}\b{\xi_0, S \xi^*_0} 
  = t^{-1}\b{S\xi_0, \xi^*_0}  
  = t^{-1}\b{\xi_0, \xi^*_0}.
\]
Using these comments line \eqref{eq:newrel} is routinely verified.
\end{proof}

\newpage

\section{Transition maps from
$\{E_i\xi^*_0\}_{i=0}^d$ and $\{E_{d-i}\xi^*_0\}_{i=0}^d$}

\indent
In Lemma \ref{lem:24bases} we gave $24$ bases for $V$. 
We are now ready to display the transition maps between ordered pairs of bases 
among these $24$.
We start with the transition maps from 
the basis $\{E_i\xi^*_0\}_{i=0}^d$ and 
the basis $\{E_{d-i}\xi^*_0\}_{i=0}^d$.

\medskip

\begin{notation}             \label{notation:main}   \samepage
Let $\Phi=(A; \{E_i\}_{i=0}^d; A^*; \{E^*_i\}_{i=0}^d)$ denote
a Leonard system in $\cal A$
 and let 
$(\{\theta_i\}_{i=0}^d; \{\theta^*_i\}_{i=0}^d;
  \{\varphi_i\}_{i=1}^d; \{\phi_i\}_{i=1}^d)$
denote the parameter array of $\Phi$.
Let $\xi_0,\xi_d,\xi^*_0,\xi^*_d$ denote nonzero vectors in $V$
that satisfy \eqref{eq:defv0vd}.
\end{notation}

\begin{theorem}             \label{thm:Eivs0}  
Referring to Notation \ref{notation:main} the following (i), (ii) hold.
\begin{itemize}
\item[(i)]
Let $\{X_i\}_{i=0}^d$ denote one of
\[
 \{E_i\}_{i=0}^d; \quad
 \{E_{d-i}\}_{i=0}^d; \quad
 \{\tau_i(A)\}_{i=0}^d; \quad
 \{\tau_{d-i}(A)\}_{i=0}^d; \quad
 \{\eta_i(A)\}_{i=0}^d; \quad
 \{\eta_{d-i}(A)\}_{i=0}^d.
 \]
Then for $0 \leq i \leq d$,
\begin{align}     
 \sum_{r=0}^d \frac{X_rE^*_0E_r}{\tr(E_rE^*_0)}
  \cdot E_i\xi^*_0 &= X_i\xi^*_0,         
                                  \label{eq:Eivs0toXivs0} \\
 \frac{\tau^*_d(\th^*_d)}{\phi}
 \frac{\b{\xi_d,\xi^*_d}}
      {\b{\xi_d,\xi^*_0}}
 \sum_{r=0}^d \tau_r(\th_r)\eta_{d-r}(\th_r)X_rE^*_dE_r
 \cdot E_i\xi^*_0 &= X_i\xi^*_d,
                                  \label{eq:Eivs0toXivsd} \\
 \sum_{r=0}^d \frac{X_{d-r}E^*_0E_r}{\tr(E_rE^*_0)}
  \cdot E_{d-i}\xi^*_0 &= X_i\xi^*_0,         
                                  \label{eq:Ed-ivs0toXivs0} \\
 \frac{\tau^*_d(\th^*_d)}{\phi}
 \frac{\b{\xi_d,\xi^*_d}}
      {\b{\xi_d,\xi^*_0}}
 \sum_{r=0}^d \tau_r(\th_r)\eta_{d-r}(\th_r)X_{d-r}E^*_dE_r
 \cdot E_{d-i}\xi^*_0 &= X_i\xi^*_d.
                                  \label{eq:Ed-ivs0toXivsd} 
\end{align}
\item[(ii)]
Let $\{X_i\}_{i=0}^d$ denote one of
\[
 \{E^*_i\}_{i=0}^d; \quad
 \{E^*_{d-i}\}_{i=0}^d; \quad
 \{\tau^*_i(A^*)\}_{i=0}^d; \quad
 \{\tau^*_{d-i}(A^*)\}_{i=0}^d; \quad
 \{\eta^*_i(A^*)\}_{i=0}^d; \quad
 \{\eta^*_{d-i}(A^*)\}_{i=0}^d.
 \]
Then for $0 \leq i \leq d$,
\begin{align}
 \frac{\b{\xi_0,\xi_0}}
      {\b{\xi_0,\xi^*_0}}
 \sum_{r=0}^d \frac{X_rE_0E^*_0E_r}{\tr(E_rE^*_0)}
 \cdot E_i\xi^*_0 &= X_i\xi_0,        \label{eq:Eivs0toXiv0}  \\
 \frac{\b{\xi_d,\xi_d}}
      {\b{\xi_d,\xi^*_0}}
 \sum_{r=0}^d \frac{X_rE_dE^*_0E_r}{\tr(E_rE^*_0)}
 \cdot E_i\xi^*_0 &= X_i\xi_d,          \label{eq:Eivs0toXivd}  \\
 \frac{\b{\xi_0,\xi_0}}
      {\b{\xi_0,\xi^*_0}}
 \sum_{r=0}^d \frac{X_{d-r}E_0E^*_0E_r}{\tr(E_rE^*_0)}
 \cdot E_{d-i}\xi^*_0 &= X_i\xi_0,        \label{eq:Ed-ivs0toXiv0}  \\
 \frac{\b{\xi_d,\xi_d}}
      {\b{\xi_d,\xi^*_0}}
 \sum_{r=0}^d \frac{X_{d-r}E_dE^*_0E_r}{\tr(E_rE^*_0)}
 \cdot E_{d-i}\xi^*_0 &= X_i\xi_d.          \label{eq:Ed-ivs0toXivd}
\end{align}
\end{itemize}
\end{theorem}

\begin{proof}
For $0 \leq i,r \leq d$ we have $E_rE_i=\delta_{r,i}E_r$.
From this and the equation on the left in \eqref{eq:EsdErEsd},
\begin{equation}         \label{eq:aux}
    \frac{E^*_0 E_r}{\tr(E_rE^*_0)}
    \cdot E_iE^*_0=\delta_{r,i} E^*_0
    \qquad\qquad (0 \leq r,i \leq d).
\end{equation}

(i):
To get \eqref{eq:Eivs0toXivs0},
multiply each side of \eqref{eq:aux} on the left by $X_r$ and on the right 
by $\xi^*_0$;
now sum the resulting equation over $r=0,\ldots,d$
and simplify using $E^*_0\xi^*_0=\xi^*_0$.
To get \eqref{eq:Eivs0toXivsd},
multiply each side of \eqref{eq:aux} on the left by $X_rE^*_dE_d$ and on 
the right by $\xi^*_0$;
now sum the resulting equation over $r=0,\ldots,d$
and simplify using
$E^*_0\xi^*_0=\xi^*_0$, 
the equations on the right in \eqref{eq:vs0}, \eqref{eq:vd},
equations \eqref{eq:trErEs0}, \eqref{eq:dds}, 
the equation on the right in \eqref{eq:trEdEs0}, 
and \eqref{eq:EsdEdEs0Er}.
To get \eqref{eq:Ed-ivs0toXivs0}, 
in line \eqref{eq:Eivs0toXivs0} replace $i$ by $d-i$
and $X_0,\ldots,X_d$ by $X_d,\ldots,X_0$.
To get \eqref{eq:Ed-ivs0toXivsd},
in line \eqref{eq:Eivs0toXivsd} replace $i$ by $d-i$
and $X_0,\ldots,X_d$ by $X_d,\ldots,X_0$.

(ii):
To get \eqref{eq:Eivs0toXiv0},
multiply each side of \eqref{eq:aux} on the left by $X_rE_0$ and on the 
right by $\xi^*_0$;
now sum the resulting equation over $r=0,\ldots,d$
and simplify using $E^*_0\xi^*_0=\xi^*_0$ and
the equation on the left in \eqref{eq:vs0}.
To get \eqref{eq:Eivs0toXivd},
multiply each side of \eqref{eq:aux} on the left by $X_rE_d$ and on the 
right by $\xi^*_0$;
now sum the resulting equation over $r=0,\ldots,d$
and simplify using $E^*_0\xi^*_0=\xi^*_0$
and the equation on the right in \eqref{eq:vs0}.
To get \eqref{eq:Ed-ivs0toXiv0},
in line \eqref{eq:Eivs0toXiv0} replace $i$ by $d-i$
and $X_0,\ldots,X_d$ by $X_d,\ldots,X_0$.
To get \eqref{eq:Ed-ivs0toXivd},
in line \eqref{eq:Eivs0toXivd}  replace $i$ by $d-i$
and $X_0,\ldots,X_d$ by $X_d,\ldots,X_0$.
\end{proof}

\newpage

\section{Transition maps from 
$\{E_i\xi^*_d\}_{i=0}^d$ and $\{E_{d-i}\xi^*_d\}_{i=0}^d$}

\begin{theorem}             \label{thm:Eivsd}   \samepage
Referring to Notation \ref{notation:main} the following (i), (ii) hold.
\begin{itemize}
\item[(i)]
Let $\{X_i\}_{i=0}^d$ denote one of
\[
 \{E_i\}_{i=0}^d; \quad \{E_{d-i}\}_{i=0}^d; \quad 
 \{\tau_i(A)\}_{i=0}^d; \quad \{\tau_{d-i}(A)\}_{i=0}^d;  \quad 
 \{\eta_i(A)\}_{i=0}^d; \quad \{\eta_{d-i}(A)\}_{i=0}^d.
\]
Then for $0 \leq i \leq d$,
\begin{align*}    
 \frac{\eta^*_d(\th^*_0)}{\vphi}
 \frac{\b{\xi_d,\xi^*_0}}
      {\b{\xi_d,\xi^*_d}}
 \sum_{r=0}^d \tau_r(\th_r)\eta_{d-r}(\th_r)X_rE^*_0E_r
 \cdot  E_i\xi^*_d &= X_i\xi^*_0,         \\
 \sum_{r=0}^d \frac{X_rE^*_dE_r}{\tr(E_rE^*_d)}
 \cdot E_i\xi^*_d &= X_i\xi^*_d,               \\
 \frac{\eta^*_d(\th^*_0)}{\vphi}
 \frac{\b{\xi_d,\xi^*_0}}
      {\b{\xi_d,\xi^*_d}}
 \sum_{r=0}^d \tau_r(\th_r)\eta_{d-r}(\th_r)X_{d-r}E^*_0E_r
 \cdot  E_{d-i}\xi^*_d &= X_i\xi^*_0,             \\
 \sum_{r=0}^d \frac{X_{d-r}E^*_dE_r}{\tr(E_rE^*_d)}
 \cdot E_{d-i}\xi^*_d &= X_i\xi^*_d.               
\end{align*}
\item[(ii)]
Let $\{X_i\}_{i=0}^d$ denote one of
\[
 \{E^*_i\}_{i=0}^d; \quad \{E^*_{d-i}\}_{i=0}^d; \quad 
 \{\tau^*_i(A^*)\}_{i=0}^d; \quad \{\tau^*_{d-i}(A^*)\}_{i=0}^d;  \quad 
 \{\eta^*_i(A^*)\}_{i=0}^d; \quad \{\eta^*_{d-i}(A^*)\}_{i=0}^d.
\]
Then for $0 \leq i \leq d$,
\begin{align*}
 \frac{\b{\xi_0,\xi_0}}
      {\b{\xi_0,\xi^*_d}}
 \sum_{r=0}^d \frac{X_rE_0E^*_dE_r}{\tr(E_rE^*_d)}
 \cdot E_i\xi^*_d  &= X_i\xi_0,              \\
 \frac{\b{\xi_d,\xi_d}}
      {\b{\xi_d,\xi^*_d}}
 \sum_{r=0}^d \frac{X_rE_dE^*_dE_r}{\tr(E_rE^*_d)}
 \cdot E_i\xi^*_d &= X_i\xi_d,               \\
 \frac{\b{\xi_0,\xi_0}}
      {\b{\xi_0,\xi^*_d}}
 \sum_{r=0}^d \frac{X_{d-r}E_0E^*_dE_r}{\tr(E_rE^*_d)}
 \cdot E_{d-i}\xi^*_d  &= X_i\xi_0,            \\
 \frac{\b{\xi_d,\xi_d}}
      {\b{\xi_d,\xi^*_d}}
 \sum_{r=0}^d \frac{X_{d-r}E_dE^*_dE_r}{\tr(E_rE^*_d)}
 \cdot E_{d-i}\xi^*_d &= X_i\xi_d.
\end{align*}
\end{itemize}
\end{theorem}

\begin{proof}
Apply Theorem \ref{thm:Eivs0} to $\Phi^\downarrow$.
\end{proof}

\section{Transition maps from 
$\{E^*_i\xi_0\}_{i=0}^d$ and $\{E^*_{d-i}\xi_0\}_{i=0}^d$}

\begin{theorem}             \label{thm:Esiv0}   \samepage
Referring to Notation \ref{notation:main} the following (i), (ii) hold.
\begin{itemize}
\item[(i)]
Let $\{X_i\}_{i=0}^d$ denote one of
\[
 \{E_i\}_{i=0}^d; \quad \{E_{d-i}\}_{i=0}^d; \quad 
 \{\tau_i(A)\}_{i=0}^d; \quad \{\tau_{d-i}(A)\}_{i=0}^d;  \quad 
 \{\eta_i(A)\}_{i=0}^d; \quad \{\eta_{d-i}(A)\}_{i=0}^d.
\]
Then for $0 \leq i \leq d$,
\begin{align*}  
 \frac{\b{\xi^*_0,\xi^*_0}}
      {\b{\xi_0,\xi^*_0}}
 \sum_{r=0}^d \frac{X_rE^*_0E_0E^*_r}{\tr(E^*_rE_0)}
 \cdot E^*_i\xi_0 &= X_i\xi^*_0,          \\
 \frac{\b{\xi^*_d,\xi^*_d}}
      {\b{\xi_0,\xi^*_d}}
 \sum_{r=0}^d \frac{X_rE^*_dE_0E^*_r}{\tr(E^*_rE_0)}
 \cdot  E^*_i\xi_0 &= X_i\xi^*_d,       \\
 \frac{\b{\xi^*_0,\xi^*_0}}
      {\b{\xi_0,\xi^*_0}}
 \sum_{r=0}^d \frac{X_{d-r}E^*_0E_0E^*_r}{\tr(E^*_rE_0)}
 \cdot E^*_{d-i}\xi_0 &= X_i\xi^*_0,         \\
 \frac{\b{\xi^*_d,\xi^*_d}}
      {\b{\xi_0,\xi^*_d}}
 \sum_{r=0}^d \frac{X_{d-r}E^*_dE_0E^*_r}{\tr(E^*_rE_0)}
 \cdot  E^*_{d-i}\xi_0 &= X_i\xi^*_d. 
\end{align*}
\item[(ii)]
Let $\{X_i\}_{i=0}^d$ denote one of
\[
 \{E^*_i\}_{i=0}^d; \quad \{E^*_{d-i}\}_{i=0}^d; \quad 
 \{\tau^*_i(A^*)\}_{i=0}^d; \quad \{\tau^*_{d-i}(A^*)\}_{i=0}^d;  \quad 
 \{\eta^*_i(A^*)\}_{i=0}^d; \quad \{\eta^*_{d-i}(A^*)\}_{i=0}^d.
\]
Then for $0 \leq i \leq d$,
\begin{align*}
 \sum_{r=0}^d \frac{X_rE_0E^*_r}{\tr(E^*_rE_0)}
 \cdot E^*_i\xi_0 &= X_i\xi_0,         \\
 \frac{\tau_d(\th_d)}{\phi}
 \frac{\b{\xi_d,\xi^*_d}}
      {\b{\xi_0,\xi^*_d}}
 \sum_{r=0}^d \tau^*_r(\th^*_r)\eta^*_{d-r}(\th^*_r)X_rE_dE^*_r
 \cdot E^*_i\xi_0 &= X_i\xi_d,          \\
 \sum_{r=0}^d \frac{X_{d-r}E_0E^*_r}{\tr(E^*_rE_0)}
 \cdot E^*_{d-i}\xi_0 &= X_i\xi_0,       \\
 \frac{\tau_d(\th_d)}{\phi}
 \frac{\b{\xi_d,\xi^*_d}}
      {\b{\xi_0,\xi^*_d}}
 \sum_{r=0}^d \tau^*_r(\th^*_r)\eta^*_{d-r}(\th^*_r)X_{d-r}E_dE^*_r
 \cdot E^*_{d-i}\xi_0 &= X_i\xi_d.      
\end{align*}
\end{itemize}
\end{theorem}

\begin{proof}
Apply Theorem \ref{thm:Eivs0} to $\Phi^*$.
\end{proof}

\section{Transition maps from 
$\{E^*_i\xi_d\}_{i=0}^d$ and $\{E^*_{d-i}\xi_d\}_{i=0}^d$}

\begin{theorem}             \label{thm:Esivd}   \samepage
Referring to Notation \ref{notation:main} the following (i), (ii) hold.
\begin{itemize}
\item[(i)]
Let $\{X_i\}_{i=0}^d$ denote one of
\[
 \{E_i\}_{i=0}^d; \quad \{E_{d-i}\}_{i=0}^d; \quad 
 \{\tau_i(A)\}_{i=0}^d; \quad \{\tau_{d-i}(A)\}_{i=0}^d;  \quad 
 \{\eta_i(A)\}_{i=0}^d; \quad \{\eta_{d-i}(A)\}_{i=0}^d.
\]
Then for $0 \leq i \leq d$,
\begin{align*}     
 \frac{\b{\xi^*_0,\xi^*_0}}
      {\b{\xi_d,\xi^*_0}}
 \sum_{r=0}^d \frac{X_rE^*_0E_dE^*_r}{\tr(E^*_rE_d)}
 \cdot E^*_i\xi_d &= X_i\xi^*_0,            \\
 \frac{\b{\xi^*_d,\xi^*_d}}
      {\b{\xi_d,\xi^*_d}}
 \sum_{r=0}^d \frac{X_rE^*_dE_dE^*_r}{\tr(E^*_rE_d)}
 \cdot E^*_i\xi_d &= X_i\xi^*_d,             \\
 \frac{\b{\xi^*_0,\xi^*_0}}
      {\b{\xi_d,\xi^*_0}}
 \sum_{r=0}^d \frac{X_{d-r}E^*_0E_dE^*_r}{\tr(E^*_rE_d)}
 \cdot E^*_{d-i}\xi_d &= X_i\xi^*_0,           \\
 \frac{\b{\xi^*_d,\xi^*_d}}
      {\b{\xi_d,\xi^*_d}}
 \sum_{r=0}^d \frac{X_{d-r}E^*_dE_dE^*_r}{\tr(E^*_rE_d)}
 \cdot E^*_{d-i}\xi_d &= X_i\xi^*_d.  
\end{align*}
\item[(ii)]
Let $\{X_i\}_{i=0}^d$ denote one of
\[
 \{E^*_i\}_{i=0}^d; \quad \{E^*_{d-i}\}_{i=0}^d; \quad 
 \{\tau^*_i(A^*)\}_{i=0}^d; \quad \{\tau^*_{d-i}(A^*)\}_{i=0}^d;  \quad 
 \{\eta^*_i(A^*)\}_{i=0}^d; \quad \{\eta^*_{d-i}(A^*)\}_{i=0}^d.
\]
Then for $0 \leq i \leq d$,
\begin{align*}
 \frac{\eta_d(\th_0)}{\vphi}
 \frac{\b{\xi_0,\xi^*_d}}
      {\b{\xi_d,\xi^*_d}}
 \sum_{r=0}^d \tau^*_r(\th^*_r)\eta^*_{d-r}(\th^*_r)X_rE_0E^*_r
 \cdot E^*_i\xi_d &= X_i\xi_0,          \\
 \sum_{r=0}^d \frac{X_rE_dE^*_r}{\tr(E^*_rE_d)}
 \cdot E^*_i\xi_d &= X_i\xi_d,             \\
 \frac{\eta_d(\th_0)}{\vphi}
 \frac{\b{\xi_0,\xi^*_d}}
      {\b{\xi_d,\xi^*_d}}
 \sum_{r=0}^d \tau^*_r(\th^*_r)\eta^*_{d-r}(\th^*_r)X_{d-r}E_0E^*_r
 \cdot E^*_{d-i}\xi_d &= X_i\xi_0,             \\
 \sum_{r=0}^d \frac{X_{d-r}E_dE^*_r}{\tr(E^*_rE_d)}
 \cdot E^*_{d-i}\xi_d &= X_i\xi_d.
\end{align*}
\end{itemize}
\end{theorem}

\begin{proof}
Apply Theorem \ref{thm:Esiv0} to $\Phi^\Downarrow$.
\end{proof}

\section{Transition maps from 
$\{\tau_i(A)\xi^*_0\}_{i=0}^d$ and $\{\tau_{d-i}(A)\xi^*_0\}_{i=0}^d$}

\begin{theorem}             \label{thm:1}   \samepage
Referring to Notation \ref{notation:main} the following (i), (ii) hold.
\begin{itemize}
\item[(i)]
Let $\{X_i\}_{i=0}^d$ denote one of
\[
 \{E_i\}_{i=0}^d; \quad \{E_{d-i}\}_{i=0}^d; \quad 
 \{\tau_i(A)\}_{i=0}^d; \quad \{\tau_{d-i}(A)\}_{i=0}^d;  \quad 
 \{\eta_i(A)\}_{i=0}^d; \quad \{\eta_{d-i}(A)\}_{i=0}^d.
\]
Then for $0 \leq i \leq d$,
\begin{align} 
 \frac{1}{\tr(E_0E^*_0)}
 \sum_{r=0}^d 
 \frac{X_rE^*_0E_0\tau^*_r(A^*)}
      {\varphi_1\varphi_2\cdots\varphi_r}
 \cdot  \tau_i(A)\xi^*_0 &= X_i\xi^*_0,       \label{eq:tauiAvs0toXivs0} \\
 \frac{\tau^*_d(\th^*_d)}{\vphi}
 \frac{\b{\xi_0,\xi^*_d}}
      {\b{\xi_0,\xi^*_0}}
 \sum_{r=0}^d X_rE^*_d\eta_{d-r}(A)
 \cdot \tau_i(A)\xi^*_0 &= X_i\xi^*_d,        \label{eq:tauiAvs0toXivsd} \\
 \frac{1}{\tr(E_0E^*_0)}
 \sum_{r=0}^d 
 \frac{X_{d-r}E^*_0E_0\tau^*_r(A^*)}
      {\varphi_1\varphi_2\cdots\varphi_r}
 \cdot  \tau_{d-i}(A)\xi^*_0 &= X_i\xi^*_0,   \label{eq:taud-iAvs0toXivs0} \\
 \frac{\tau^*_d(\th^*_d)}{\vphi}
 \frac{\b{\xi_0,\xi^*_d}}
      {\b{\xi_0,\xi^*_0}}
 \sum_{r=0}^d X_{d-r}E^*_d\eta_{d-r}(A)
 \cdot \tau_{d-i}(A)\xi^*_0 &= X_i\xi^*_d.        \label{eq:taud-iAvs0toXivsd} 
\end{align}
\item[(ii)]
Let $\{X_i\}_{i=0}^d$ denote one of
\[
 \{E^*_i\}_{i=0}^d; \quad \{E^*_{d-i}\}_{i=0}^d; \quad 
 \{\tau^*_i(A^*)\}_{i=0}^d; \quad \{\tau^*_{d-i}(A^*)\}_{i=0}^d;  \quad 
 \{\eta^*_i(A^*)\}_{i=0}^d; \quad \{\eta^*_{d-i}(A^*)\}_{i=0}^d.
\]
Then for $0 \leq i \leq d$,
\begin{align}
 \frac{\b{\xi_0,\xi_0}}
      {\b{\xi_0,\xi^*_0}}
 \sum_{r=0}^d 
 \frac{X_rE_0\tau^*_r(A^*)}{\varphi_1\varphi_2\cdots\varphi_r}
 \cdot \tau_i(A)\xi^*_0 &= X_i\xi_0,          \label{eq:tauiAvs0toXiv0} \\
 \frac{\tau^*_d(\th^*_d)}{\phi}
 \frac{\b{\xi_d,\xi_d}}
      {\b{\xi_d,\xi^*_0}}
 \sum_{r=0}^d X_rE_dE^*_d\eta_{d-r}(A)
 \cdot \tau_i(A)\xi^*_0 &= X_i\xi_d,         \label{eq:tauiAvs0toXivd}  \\
 \frac{\b{\xi_0,\xi_0}}
      {\b{\xi_0,\xi^*_0}}
 \sum_{r=0}^d 
 \frac{X_{d-r}E_0\tau^*_r(A^*)}{\varphi_1\varphi_2\cdots\varphi_r}
 \cdot \tau_{d-i}(A)\xi^*_0 &= X_i\xi_0,          \label{eq:taud-iAvs0toXiv0} \\
 \frac{\tau^*_d(\th^*_d)}{\phi}
 \frac{\b{\xi_d,\xi_d}}
      {\b{\xi_d,\xi^*_0}}
 \sum_{r=0}^d X_{d-r}E_dE^*_d\eta_{d-r}(A)
 \cdot \tau_{d-i}(A)\xi^*_0 &= X_i\xi_d.         \label{eq:taud-iAvs0toXivd}
\end{align}
\end{itemize}
\end{theorem}

\begin{proof}
Multiply both sides of \eqref{eq:basic} on the left by $E^*_0$, 
use the equation on the left in \eqref{eq:EsdErEsd},
and replace $(i,j)$ with $(r,i)$ to obtain
\begin{equation}                  \label{eq:auxB}
  \frac{1}{\tr(E_0E^*_0)}
  \frac{E^*_0E_0\tau^*_r(A^*)}{\vphi_1\vphi_2\cdots\vphi_r}
  \cdot \tau_i(A)E^*_0 = \delta_{r,i} E^*_0
  \qquad\qquad (0 \leq r,i \leq d).
\end{equation}

(i):
To get \eqref{eq:tauiAvs0toXivs0},
multiply each side of \eqref{eq:auxB} on the left by $X_r$ and on the 
right by $\xi^*_0$;
now sum the resulting equation over $r=0,\ldots,d$ and
simplify using $E^*_0\xi^*_0=\xi^*_0$.
Next we show \eqref{eq:tauiAvs0toXivsd}.
By \cite[Corollary 5.3]{NT:mu} we have
\begin{equation}           \label{eq:mu}
  E^*_dE_0\tau^*_r(A^*)
  = \frac{\vphi_1\vphi_2\cdots\vphi_r}
         {\eta_d(\th_0)}
    E^*_d\eta_{d-r}(A)
\end{equation}
for $0 \leq r \leq d$.
Now to get \eqref{eq:tauiAvs0toXivsd},  
multiply each side of \eqref{eq:auxB} on the left by $X_rE^*_dE_0$ and 
on the right by $\xi^*_0$,
sum the resulting equation over $r=0,\ldots,d$ and
simplify using $E^*_0\xi^*_0=\xi^*_0$, 
the equations on the left in \eqref{eq:E0EsrE0}, \eqref{eq:vs0}, 
the equation on the right in \eqref{eq:v0}, 
equations \eqref{eq:0ds}, \eqref{eq:mu}, 
and the equation on the right in \eqref{eq:trE0Es0}.
To get \eqref{eq:taud-iAvs0toXivs0}, 
in line \eqref{eq:tauiAvs0toXivs0} replace $i$ by $d-i$
and $X_0,\ldots,X_d$ by $X_d,\ldots,X_0$.
To get  \eqref{eq:taud-iAvs0toXivsd},
in line \eqref{eq:tauiAvs0toXivsd} replace $i$ by $d-i$
and $X_0,\ldots,X_d$ by $X_d,\ldots,X_0$.

(ii):
To get \eqref{eq:tauiAvs0toXiv0}, 
multiply each side of \eqref{eq:auxB} on the left by $X_rE_0$ and on the 
right by $\xi^*_0$;
now sum the resulting equation over $r=0,\ldots,d$ and
simplify using $E^*_0\xi^*_0=\xi^*_0$ and
the equations on the left in \eqref{eq:E0EsrE0}, \eqref{eq:vs0}.
To get \eqref{eq:tauiAvs0toXivd},
multiply each side of \eqref{eq:auxB} on the left by $X_rE_dE^*_dE_0$ and 
on the right by $\xi^*_0$;
now sum the resulting equation over $r=0,\ldots,d$
and simplify using $E^*_0\xi^*_0=\xi^*_0$,
the equations on the left in \eqref{eq:E0EsrE0}, \eqref{eq:vs0}, 
the equations on the right in \eqref{eq:vsd}, \eqref{eq:v0},
equations \eqref{eq:0ds}, \eqref{eq:newrel}, \eqref{eq:mu}, 
and the equation on the right in \eqref{eq:trE0Es0}.
To get \eqref{eq:taud-iAvs0toXiv0},
in line \eqref{eq:tauiAvs0toXiv0} replace $i$ by $d-i$
and $X_0,\ldots,X_d$ by $X_d,\ldots,X_0$.
To get \eqref{eq:taud-iAvs0toXivd},
in line \eqref{eq:tauiAvs0toXivd} replace  $i$ by $d-i$
and $X_0,\ldots,X_d$ by $X_d,\ldots,X_0$.
\end{proof}

\section{Transition maps from 
$\{\eta_i(A)\xi^*_0\}_{i=0}^d$ and $\{\eta_{d-i}(A)\xi^*_0\}_{i=0}^d$}

\begin{theorem}             \label{thm:D}   \samepage
Referring to Notation \ref{notation:main} the following (i), (ii) hold.
\begin{itemize}
\item[(i)]
Let $\{X_i\}_{i=0}^d$ denote one of
\[
 \{E_i\}_{i=0}^d; \quad \{E_{d-i}\}_{i=0}^d; \quad 
 \{\tau_i(A)\}_{i=0}^d; \quad \{\tau_{d-i}(A)\}_{i=0}^d;  \quad 
 \{\eta_i(A)\}_{i=0}^d; \quad \{\eta_{d-i}(A)\}_{i=0}^d.
\]
Then for $0 \leq i \leq d$,
\begin{align*} 
 \frac{1}{\tr(E_dE^*_0)} 
 \sum_{r=0}^d
 \frac{X_rE^*_0E_d\tau^*_r(A^*)}{\phi_1\phi_2\cdots\phi_r}
 \cdot \eta_i(A)\xi^*_0 &= X_i\xi^*_0,        \\
 \frac{\tau^*_d(\th^*_d)}{\phi}
 \frac{\b{\xi_d,\xi^*_d}}
      {\b{\xi_d,\xi^*_0}}
 \sum_{r=0}^d X_rE^*_d\tau_{d-r}(A)
 \cdot \eta_i(A)\xi^*_0 &= X_i\xi^*_d,      \\
 \frac{1}{\tr(E_dE^*_0)} 
 \sum_{r=0}^d
 \frac{X_{d-r}E^*_0E_d\tau^*_r(A^*)}{\phi_1\phi_2\cdots\phi_r}
 \cdot \eta_{d-i}(A)\xi^*_0 &= X_i\xi^*_0,      \\
 \frac{\tau^*_d(\th^*_d)}{\phi}
 \frac{\b{\xi_d,\xi^*_d}}
      {\b{\xi_d,\xi^*_0}}
 \sum_{r=0}^d X_{d-r}E^*_d\tau_{d-r}(A)
 \cdot \eta_{d-i}(A)\xi^*_0 &= X_i\xi^*_d. 
\end{align*}
\item[(ii)]
Let $\{X_i\}_{i=0}^d$ denote one of
\[
 \{E^*_i\}_{i=0}^d; \quad \{E^*_{d-i}\}_{i=0}^d; \quad 
 \{\tau^*_i(A^*)\}_{i=0}^d; \quad \{\tau^*_{d-i}(A^*)\}_{i=0}^d;  \quad 
 \{\eta^*_i(A^*)\}_{i=0}^d; \quad \{\eta^*_{d-i}(A^*)\}_{i=0}^d.
\]
Then for $0 \leq i \leq d$,
\begin{align*}
 \frac{\tau^*_d(\th^*_d)}{\vphi}
 \frac{\b{\xi_0,\xi_0}}
      {\b{\xi_0,\xi^*_0}}
 \sum_{r=0}^d X_rE_0E^*_d\tau_{d-r}(A)
 \cdot \eta_i(A)\xi^*_0 &= X_i\xi_0,        \\
 \frac{\b{\xi_d,\xi_d}}
      {\b{\xi_d,\xi^*_0}}
 \sum_{r=0}^d 
 \frac{X_rE_d\tau^*_r(A^*)}{\phi_1\phi_2\cdots\phi_r}
 \cdot  \eta_i(A)\xi^*_0 &= X_i\xi_d,       \\
 \frac{\tau^*_d(\th^*_d)}{\vphi}
 \frac{\b{\xi_0,\xi_0}}
      {\b{\xi_0,\xi^*_0}}
 \sum_{r=0}^d X_{d-r}E_0E^*_d\tau_{d-r}(A)
 \cdot \eta_{d-i}(A)\xi^*_0 &= X_i\xi_0,    \\
 \frac{\b{\xi_d,\xi_d}}
      {\b{\xi_d,\xi^*_0}}
 \sum_{r=0}^d 
 \frac{X_{d-r}E_d\tau^*_r(A^*)}{\phi_1\phi_2\cdots\phi_r}
 \cdot  \eta_{d-i}(A)\xi^*_0 &= X_i\xi_d. 
\end{align*}
\end{itemize}
\end{theorem}

\begin{proof}
Apply Theorem \ref{thm:1} to $\Phi^\Downarrow$.
\end{proof}

\section{Transition maps from 
$\{\tau_i(A)\xi^*_d\}_{i=0}^d$ and  
$\{\tau_{d-i}(A)\xi^*_d\}_{i=0}^d$}

\begin{theorem}             \label{thm:d}   \samepage
Referring to Notation \ref{notation:main} the following (i), (ii) hold.
\begin{itemize}
\item[(i)]
Let $\{X_i\}_{i=0}^d$ denote one of
\[
 \{E_i\}_{i=0}^d; \quad \{E_{d-i}\}_{i=0}^d; \quad 
 \{\tau_i(A)\}_{i=0}^d; \quad \{\tau_{d-i}(A)\}_{i=0}^d;  \quad 
 \{\eta_i(A)\}_{i=0}^d; \quad \{\eta_{d-i}(A)\}_{i=0}^d.
\]
Then for $0 \leq i \leq d$,
\begin{align*}
 \frac{\eta^*_d(\th^*_0)}{\phi}
 \frac{\b{\xi_0,\xi^*_0}}
      {\b{\xi_0,\xi^*_d}}
 \sum_{r=0}^d X_rE^*_0\eta_{d-r}(A)
 \cdot  \tau_i(A)\xi^*_d &= X_i\xi^*_0,       \\
 \frac{1}{\tr(E_0E^*_d)}
 \sum_{r=0}^d 
 \frac{X_rE^*_dE_0\eta^*_{r}(A^*)}
      {\phi_d\phi_{d-1}\cdots\phi_{d-r+1}}
 \cdot  \tau_i(A)\xi^*_d &=  X_i\xi^*_d,    \\
 \frac{\eta^*_d(\th^*_0)}{\phi}
 \frac{\b{\xi_0,\xi^*_0}}
      {\b{\xi_0,\xi^*_d}}
 \sum_{r=0}^d X_{d-r}E^*_0\eta_{d-r}(A)
 \cdot  \tau_{d-i}(A)\xi^*_d &= X_i\xi^*_0,     \\
 \frac{1}{\tr(E_0E^*_d)}
 \sum_{r=0}^d 
 \frac{X_{d-r}E^*_dE_0\eta^*_{r}(A^*)}
      {\phi_d\phi_{d-1}\cdots\phi_{d-r+1}}
 \cdot  \tau_{d-i}(A)\xi^*_d &=  X_i\xi^*_d. 
\end{align*}
\item[(ii)]
Let $\{X_i\}_{i=0}^d$ denote one of
\[
 \{E^*_i\}_{i=0}^d; \quad \{E^*_{d-i}\}_{i=0}^d; \quad 
 \{\tau^*_i(A^*)\}_{i=0}^d; \quad \{\tau^*_{d-i}(A^*)\}_{i=0}^d;  \quad 
 \{\eta^*_i(A^*)\}_{i=0}^d; \quad \{\eta^*_{d-i}(A^*)\}_{i=0}^d.
\]
Then for $0 \leq i \leq d$,
\begin{align*}
 \frac{\b{\xi_0,\xi_0}}
      {\b{\xi_0,\xi^*_d}}
 \sum_{r=0}^d 
  \frac{X_rE_0\eta^*_r(A^*)}{\phi_{d}\phi_{d-1}\cdots\phi_{d-r+1}}
 \cdot \tau_i(A)\xi^*_d &= X_i\xi_0,           \\
 \frac{\eta^*_d(\th^*_0)}{\vphi}
 \frac{\b{\xi_d,\xi_d}}
       {\b{\xi_d,\xi^*_d}}
 \sum_{r=0}^d X_rE_dE^*_0\eta_{d-r}(A)
 \cdot \tau_i(A)\xi^*_d &= X_i\xi_d,      \\
 \frac{\b{\xi_0,\xi_0}}
      {\b{\xi_0,\xi^*_d}}
 \sum_{r=0}^d 
  \frac{X_{d-r}E_0\eta^*_r(A^*)}{\phi_{d}\phi_{d-1}\cdots\phi_{d-r+1}}
 \cdot \tau_{d-i}(A)\xi^*_d &= X_i\xi_0,       \\
 \frac{\eta^*_d(\th^*_0)}{\vphi}
 \frac{\b{\xi_d,\xi_d}}
       {\b{\xi_d,\xi^*_d}}
 \sum_{r=0}^d X_{d-r}E_dE^*_0\eta_{d-r}(A)
 \cdot \tau_{d-i}(A)\xi^*_d &= X_i\xi_d.    
\end{align*}
\end{itemize}
\end{theorem}

\begin{proof}
Apply Theorem \ref{thm:1} to $\Phi^\downarrow$.
\end{proof}

\section{Transition maps from 
$\{\eta_i(A)\xi^*_d\}_{i=0}^d$ and
$\{\eta_{d-i}(A)\xi^*_d\}_{i=0}^d$}

\begin{theorem}             \label{thm:dD}   \samepage
Referring to Notation \ref{notation:main} the following (i), (ii) hold.
\begin{itemize}
\item[(i)]
Let $\{X_i\}_{i=0}^d$ denote one of
\[
 \{E_i\}_{i=0}^d; \quad \{E_{d-i}\}_{i=0}^d; \quad 
 \{\tau_i(A)\}_{i=0}^d; \quad \{\tau_{d-i}(A)\}_{i=0}^d;  \quad 
 \{\eta_i(A)\}_{i=0}^d; \quad \{\eta_{d-i}(A)\}_{i=0}^d.
\]
Then for $0 \leq i \leq d$,
\begin{align*} 
 \frac{\eta^*_d(\th^*_0)}{\vphi}
 \frac{\b{\xi_d,\xi^*_0}}
      {\b{\xi_d,\xi^*_d}}
 \sum_{r=0}^d X_rE^*_0\tau_{d-r}(A)
 \cdot \eta_i(A)\xi^*_d &= X_i\xi^*_0,       \\
 \frac{1}{\tr(E_dE^*_d)}
 \sum_{r=0}^d 
 \frac{X_rE^*_dE_d\eta^*_{r}(A^*)}
      {\varphi_d\varphi_{d-1}\cdots\varphi_{d-r+1}}
 \cdot \eta_i(A)\xi^*_d &= X_i\xi^*_d,          \\
 \frac{\eta^*_d(\th^*_0)}{\vphi}
 \frac{\b{\xi_d,\xi^*_0}}
      {\b{\xi_d,\xi^*_d}}
 \sum_{r=0}^d X_{d-r}E^*_0\tau_{d-r}(A)
 \cdot \eta_{d-i}(A)\xi^*_d &= X_i\xi^*_0,   \\
 \frac{1}{\tr(E_dE^*_d)}
 \sum_{r=0}^d 
 \frac{X_{d-r}E^*_dE_d\eta^*_{r}(A^*)}
      {\varphi_d\varphi_{d-1}\cdots\varphi_{d-r+1}}
 \cdot \eta_{d-i}(A)\xi^*_d &= X_i\xi^*_d.  
\end{align*}
\item[(ii)]
Let $\{X_i\}_{i=0}^d$ denote one of
\[
 \{E^*_i\}_{i=0}^d; \quad \{E^*_{d-i}\}_{i=0}^d; \quad 
 \{\tau^*_i(A^*)\}_{i=0}^d; \quad \{\tau^*_{d-i}(A^*)\}_{i=0}^d;  \quad 
 \{\eta^*_i(A^*)\}_{i=0}^d; \quad \{\eta^*_{d-i}(A^*)\}_{i=0}^d.
\]
Then for $0 \leq i \leq d$,
\begin{align*}
 \frac{\eta^*_d(\th^*_0)}{\phi}
 \frac{\b{\xi_0,\xi_0}}
      {\b{\xi_0,\xi^*_d}}
 \sum_{r=0}^d X_rE_0E^*_0\tau_{d-r}(A)
 \cdot \eta_i(A)\xi^*_d &= X_i\xi_0,     \\
 \frac{\b{\xi_d,\xi_d}}
      {\b{\xi_d,\xi^*_d}}
 \sum_{r=0}^d 
  \frac{X_rE_d\eta^*_{r}(A^*)}{\varphi_d\varphi_{d-1}\cdots\varphi_{d-r+1}}
 \cdot \eta_i(A)\xi^*_d &= X_i\xi_d,         \\
 \frac{\eta^*_d(\th^*_0)}{\phi}
 \frac{\b{\xi_0,\xi_0}}
      {\b{\xi_0,\xi^*_d}}
 \sum_{r=0}^d X_{d-r}E_0E^*_0\tau_{d-r}(A)
 \cdot \eta_{d-i}(A)\xi^*_d &= X_i\xi_0,    \\
 \frac{\b{\xi_d,\xi_d}}
      {\b{\xi_d,\xi^*_d}}
 \sum_{r=0}^d 
  \frac{X_{d-r}E_d\eta^*_{r}(A^*)}{\varphi_d\varphi_{d-1}\cdots\varphi_{d-r+1}}
 \cdot \eta_{d-i}(A)\xi^*_d &= X_i\xi_d.      
\end{align*}
\end{itemize}
\end{theorem}

\begin{proof}
Apply Theorem \ref{thm:1} to $\Phi^{\downarrow\Downarrow}$.
\end{proof}

\section{Transition maps from 
$\{\tau^*_i(A^*)\xi_0\}_{i=0}^d$ and 
$\{\tau^*_{d-i}(A^*)\xi_0\}_{i=0}^d$}

\begin{theorem}             \label{thm:s}   \samepage
Referring to Notation \ref{notation:main} the following (i), (ii) hold.
\begin{itemize}
\item[(i)]
Let $\{X_i\}_{i=0}^d$ denote one of
\[
 \{E_i\}_{i=0}^d; \quad \{E_{d-i}\}_{i=0}^d; \quad 
 \{\tau_i(A)\}_{i=0}^d; \quad \{\tau_{d-i}(A)\}_{i=0}^d;  \quad 
 \{\eta_i(A)\}_{i=0}^d; \quad \{\eta_{d-i}(A)\}_{i=0}^d.
\]
Then for $0 \leq i \leq d$,
\begin{align*} 
 \frac{\b{\xi^*_0,\xi^*_0}}
      {\b{\xi_0,\xi^*_0}}
 \sum_{r=0}^d 
   \frac{X_rE^*_0\tau_{r}(A)}{\varphi_1\varphi_2\cdots\varphi_r}
 \cdot \tau^*_i(A^*)\xi_0 &= X_i\xi^*_0,      \\
 \frac{\tau_d(\th_d)}{\phi}
 \frac{\b{\xi^*_d,\xi^*_d}}
      {\b{\xi_0,\xi^*_d}}
 \sum_{r=0}^d X_rE^*_dE_d\eta^*_{d-r}(A^*)
 \cdot \tau^*_i(A^*)\xi_0 &= X_i\xi^*_d,     \\
 \frac{\b{\xi^*_0,\xi^*_0}}
      {\b{\xi_0,\xi^*_0}}
 \sum_{r=0}^d 
   \frac{X_{d-r}E^*_0\tau_{r}(A)}{\varphi_1\varphi_2\cdots\varphi_r}
 \cdot \tau^*_{d-i}(A^*)\xi_0 &= X_i\xi^*_0,   \\
 \frac{\tau_d(\th_d)}{\phi}
 \frac{\b{\xi^*_d,\xi^*_d}}
      {\b{\xi_0,\xi^*_d}}
 \sum_{r=0}^d X_{d-r}E^*_dE_d\eta^*_{d-r}(A^*)
 \cdot \tau^*_{d-i}(A^*)\xi_0 &= X_i\xi^*_d.  
\end{align*}
\item[(ii)]
Let $\{X_i\}_{i=0}^d$ denote one of
\[
 \{E^*_i\}_{i=0}^d; \quad \{E^*_{d-i}\}_{i=0}^d; \quad 
 \{\tau^*_i(A^*)\}_{i=0}^d; \quad \{\tau^*_{d-i}(A^*)\}_{i=0}^d;  \quad 
 \{\eta^*_i(A^*)\}_{i=0}^d; \quad \{\eta^*_{d-i}(A^*)\}_{i=0}^d.
\]
Then for $0 \leq i \leq d$,
\begin{align*}
 \frac{1}{\tr(E_0E^*_0)}
 \sum_{r=0}^d
 \frac{X_rE_0E^*_0\tau_r(A)}
      {\varphi_1\varphi_2\cdots\varphi_r}
 \cdot \tau^*_i(A^*)\xi_0 &= X_i\xi_0,      \\
 \frac{\tau_d(\th_d)}{\vphi}
 \frac{\b{\xi_d,\xi^*_0}}
      {\b{\xi_0,\xi^*_0}}
 \sum_{r=0}^d X_rE_d\eta^*_{d-r}(A^*)
 \cdot \tau^*_i(A^*)\xi_0 &= X_i\xi_d,      \\
 \frac{1}{\tr(E_0E^*_0)}
 \sum_{r=0}^d
 \frac{X_{d-r}E_0E^*_0\tau_r(A)}
      {\varphi_1\varphi_2\cdots\varphi_r}
 \cdot \tau^*_{d-i}(A^*)\xi_0 &= X_i\xi_0,      \\
 \frac{\tau_d(\th_d)}{\vphi}
 \frac{\b{\xi_d,\xi^*_0}}
      {\b{\xi_0,\xi^*_0}}
 \sum_{r=0}^d X_{d-r}E_d\eta^*_{d-r}(A^*)
 \cdot \tau^*_{d-i}(A^*)\xi_0 &= X_i\xi_d. 
\end{align*}
\end{itemize}
\end{theorem}

\begin{proof}
Apply Theorem \ref{thm:1} to $\Phi^*$.
\end{proof}

\section{Transition maps from 
$\{\eta^*_i(A^*)\xi_0\}_{i=0}^d$ and 
$\{\eta^*_{d-i}(A^*)\xi_0\}_{i=0}^d$}

\begin{theorem}             \label{thm:Ds}   \samepage
Referring to Notation \ref{notation:main} the following (i), (ii) hold.
\begin{itemize}
\item[(i)]
Let $\{X_i\}_{i=0}^d$ denote one of
\[
 \{E_i\}_{i=0}^d; \quad \{E_{d-i}\}_{i=0}^d; \quad 
 \{\tau_i(A)\}_{i=0}^d; \quad \{\tau_{d-i}(A)\}_{i=0}^d;  \quad 
 \{\eta_i(A)\}_{i=0}^d; \quad \{\eta_{d-i}(A)\}_{i=0}^d.
\]
Then for $0 \leq i \leq d$,
\begin{align*}
 \frac{\tau_d(\th_d)}{\vphi}
 \frac{\b{\xi^*_0,\xi^*_0}}
      {\b{\xi_0,\xi^*_0}}
 \sum_{r=0}^d X_rE^*_0E_d\tau^*_{d-r}(A^*)
 \cdot \eta^*_i(A^*)\xi_0 &= X_i\xi^*_0,     \\
 \frac{\b{\xi^*_d,\xi^*_d}}
      {\b{\xi_0,\xi^*_d}}
 \sum_{r=0}^d 
  \frac{X_rE^*_d\tau_{r}(A)}{\phi_d\phi_{d-1}\cdots\phi_{d-r+1}}
 \cdot \eta^*_i(A^*)\xi_0 &= X_i\xi^*_d,    \\
 \frac{\tau_d(\th_d)}{\vphi}
 \frac{\b{\xi^*_0,\xi^*_0}}
      {\b{\xi_0,\xi^*_0}}
 \sum_{r=0}^d X_{d-r}E^*_0E_d\tau^*_{d-r}(A^*)
 \cdot \eta^*_{d-i}(A^*)\xi_0 &= X_i\xi^*_0,   \\
 \frac{\b{\xi^*_d,\xi^*_d}}
      {\b{\xi_0,\xi^*_d}}
 \sum_{r=0}^d 
  \frac{X_{d-r}E^*_d\tau_{r}(A)}{\phi_d\phi_{d-1}\cdots\phi_{d-r+1}}
 \cdot \eta^*_{d-i}(A^*)\xi_0 &= X_i\xi^*_d.   
\end{align*}
\item[(ii)]
Let $\{X_i\}_{i=0}^d$ denote one of
\[
 \{E^*_i\}_{i=0}^d; \quad \{E^*_{d-i}\}_{i=0}^d; \quad 
 \{\tau^*_i(A^*)\}_{i=0}^d; \quad \{\tau^*_{d-i}(A^*)\}_{i=0}^d;  \quad 
 \{\eta^*_i(A^*)\}_{i=0}^d; \quad \{\eta^*_{d-i}(A^*)\}_{i=0}^d.
\]
Then for $0 \leq i \leq d$,
\begin{align*}
 \frac{1}{\tr(E_0E^*_d)}
 \sum_{r=0}^d 
  \frac{X_rE_0E^*_d\tau_{r}(A)}
       {\phi_d\phi_{d-1}\cdots\phi_{d-r+1}}
 \cdot \eta^*_i(A^*)\xi_0 &= X_i\xi_0,        \\
 \frac{\tau_d(\th_d)}{\phi}
 \frac{\b{\xi_d,\xi^*_d}}
      {\b{\xi_0,\xi^*_d}}
 \sum_{r=0}^d X_rE_d\tau^*_{d-r}(A^*)
 \cdot \eta^*_i(A^*)\xi_0 &= X_i\xi_d,      \\
 \frac{1}{\tr(E_0E^*_d)}
 \sum_{r=0}^d 
  \frac{X_{d-r}E_0E^*_d\tau_{r}(A)}
       {\phi_d\phi_{d-1}\cdots\phi_{d-r+1}}
 \cdot \eta^*_{d-i}(A^*)\xi_0 &= X_i\xi_0,    \\
 \frac{\tau_d(\th_d)}{\phi}
 \frac{\b{\xi_d,\xi^*_d}}
      {\b{\xi_0,\xi^*_d}}
 \sum_{r=0}^d  X_{d-r}E_d\tau^*_{d-r}(A^*)
 \cdot \eta^*_{d-i}(A^*)\xi_0 &= X_i\xi_d. 
\end{align*}
\end{itemize}
\end{theorem}

\begin{proof}
Apply Theorem \ref{thm:s} to $\Phi^\downarrow$.
\end{proof}

\section{Transition maps from 
$\{\tau^*_i(A^*)\xi_d\}_{i=0}^d$ and 
$\{\tau^*_{d-i}(A^*)\xi_d\}_{i=0}^d$}

\begin{theorem}             \label{thm:ds}   \samepage
Referring to Notation \ref{notation:main} the following (i), (ii) hold.
\begin{itemize}
\item[(i)]
Let $\{X_i\}_{i=0}^d$ denote one of
\[
 \{E_i\}_{i=0}^d; \quad \{E_{d-i}\}_{i=0}^d; \quad 
 \{\tau_i(A)\}_{i=0}^d; \quad \{\tau_{d-i}(A)\}_{i=0}^d;  \quad 
 \{\eta_i(A)\}_{i=0}^d; \quad \{\eta_{d-i}(A)\}_{i=0}^d.
\]
Then for $0 \leq i \leq d$,
\begin{align*} 
 \frac{\b{\xi^*_0,\xi^*_0}}
      {\b{\xi_d,\xi^*_0}}
 \sum_{r=0}^d 
   \frac{X_rE^*_0\eta_{r}(A)}{\phi_1\phi_2\cdots\phi_r}
 \cdot \tau^*_i(A^*)\xi_d &= X_i\xi^*_0,        \\
 \frac{\eta_d(\th_0)}{\vphi}
 \frac{\b{\xi^*_d,\xi^*_d}}
      {\b{\xi_d,\xi^*_d}}
 \sum_{r=0}^d X_rE^*_dE_0\eta^*_{d-r}(A^*)
 \cdot \tau^*_i(A^*)\xi_d &= X_i\xi^*_d,   \\
 \frac{\b{\xi^*_0,\xi^*_0}}
      {\b{\xi_d,\xi^*_0}}
 \sum_{r=0}^d 
   \frac{X_{d-r}E^*_0\eta_{r}(A)}{\phi_1\phi_2\cdots\phi_r}
 \cdot \tau^*_{d-i}(A^*)\xi_d &= X_i\xi^*_0,   \\
 \frac{\eta_d(\th_0)}{\vphi}
 \frac{\b{\xi^*_d,\xi^*_d}}
      {\b{\xi_d,\xi^*_d}}
 \sum_{r=0}^d X_{d-r}E^*_dE_0\eta^*_{d-r}(A^*)
 \cdot \tau^*_{d-i}(A^*)\xi_d &= X_i\xi^*_d.  
\end{align*}
\item[(ii)]
Let $\{X_i\}_{i=0}^d$ denote one of
\[
 \{E^*_i\}_{i=0}^d; \quad \{E^*_{d-i}\}_{i=0}^d; \quad 
 \{\tau^*_i(A^*)\}_{i=0}^d; \quad \{\tau^*_{d-i}(A^*)\}_{i=0}^d;  \quad 
 \{\eta^*_i(A^*)\}_{i=0}^d; \quad \{\eta^*_{d-i}(A^*)\}_{i=0}^d.
\]
Then for $0 \leq i \leq d$,
\begin{align*}
 \frac{\eta_d(\th_0)}{\phi}
 \frac{\b{\xi_0,\xi^*_0}}
      {\b{\xi_d,\xi^*_0}}
 \sum_{r=0}^d X_rE_0\eta^*_{d-r}(A^*)
 \cdot \tau^*_i(A^*)\xi_d &= X_i\xi_0,        \\
 \frac{1}{\tr(E_dE^*_0)}
 \sum_{r=0}^d 
 \frac{X_rE_dE^*_0\eta_{r}(A)}
      {\phi_1\phi_2\cdots\phi_r}
 \cdot \tau^*_i(A^*)\xi_d &= X_i\xi_d,     \\
 \frac{\eta_d(\th_0)}{\phi}
 \frac{\b{\xi_0,\xi^*_0}}
      {\b{\xi_d,\xi^*_0}}
 \sum_{r=0}^d X_{d-r}E_0\eta^*_{d-r}(A^*)
 \cdot \tau^*_{d-i}(A^*)\xi_d &= X_i\xi_0,     \\
 \frac{1}{\tr(E_dE^*_0)}
 \sum_{r=0}^d 
 \frac{X_{d-r}E_dE^*_0\eta_{r}(A)}
      {\phi_1\phi_2\cdots\phi_r}
 \cdot \tau^*_{d-i}(A^*)\xi_d &= X_i\xi_d. 
\end{align*}
\end{itemize}
\end{theorem}

\begin{proof}
Apply Theorem \ref{thm:s} to $\Phi^\Downarrow$.
\end{proof}

\section{Transition maps from 
$\{\eta^*_i(A^*)\xi_d\}_{i=0}^d$ and
$\{\eta^*_{d-i}(A^*)\xi_d\}_{i=0}^d$}

\begin{theorem}             \label{thm:dDs}   \samepage
Referring to Notation \ref{notation:main} the following (i), (ii) hold.
\begin{itemize}
\item[(i)]
Let $\{X_i\}_{i=0}^d$ denote one of
\[
 \{E_i\}_{i=0}^d; \quad \{E_{d-i}\}_{i=0}^d; \quad 
 \{\tau_i(A)\}_{i=0}^d; \quad \{\tau_{d-i}(A)\}_{i=0}^d;  \quad 
 \{\eta_i(A)\}_{i=0}^d; \quad \{\eta_{d-i}(A)\}_{i=0}^d.
\]
Then for $0 \leq i \leq d$,
\begin{align*}
 \frac{\eta_d(\th_0)}{\phi}
 \frac{\b{\xi^*_0,\xi^*_0}}
      {\b{\xi_d,\xi^*_0}}
 \sum_{r=0}^d X_rE^*_0E_0\tau^*_{d-r}(A^*)
 \cdot \eta^*_i(A^*)\xi_d &= X_i\xi^*_0,      \\
 \frac{\b{\xi^*_d,\xi^*_d}}
      {\b{\xi_d,\xi^*_d}}
 \sum_{r=0}^d 
 \frac{X_rE^*_d\eta_{r}(A)}
      {\varphi_d\varphi_{d-1}\cdots\varphi_{d-r+1}}
 \cdot \eta^*_i(A^*)\xi_d &=  X_i\xi^*_d,       \\
 \frac{\eta_d(\th_0)}{\phi}
 \frac{\b{\xi^*_0,\xi^*_0}}
      {\b{\xi_d,\xi^*_0}}
 \sum_{r=0}^d X_{d-r}E^*_0E_0\tau^*_{d-r}(A^*)
 \cdot \eta^*_{d-i}(A^*)\xi_d &= X_i\xi^*_0,   \\
 \frac{\b{\xi^*_d,\xi^*_d}}
      {\b{\xi_d,\xi^*_d}}
 \sum_{r=0}^d 
 \frac{X_{d-r}E^*_d\eta_{r}(A)}
      {\varphi_d\varphi_{d-1}\cdots\varphi_{d-r+1}}
 \cdot \eta^*_{d-i}(A^*)\xi_d &=  X_i\xi^*_d. 
\end{align*}
\item[(ii)]
Let $\{X_i\}_{i=0}^d$ denote one of
\[
 \{E^*_i\}_{i=0}^d; \quad \{E^*_{d-i}\}_{i=0}^d; \quad 
 \{\tau^*_i(A^*)\}_{i=0}^d; \quad \{\tau^*_{d-i}(A^*)\}_{i=0}^d;  \quad 
 \{\eta^*_i(A^*)\}_{i=0}^d; \quad \{\eta^*_{d-i}(A^*)\}_{i=0}^d.
\]
Then for $0 \leq i \leq d$,
\begin{align*}
 \frac{\eta_d(\th_0)}{\vphi}
 \frac{\b{\xi_0,\xi^*_d}}
      {\b{\xi_d,\xi^*_d}}
 \sum_{r=0}^d X_rE_0\tau^*_{d-r}(A^*)
 \cdot \eta^*_i(A^*)\xi_d &= X_i\xi_0,        \\
 \frac{1}{\tr(E_dE^*_d)}
 \sum_{r=0}^d 
 \frac{X_rE_dE^*_d\eta_{r}(A)}
      {\varphi_d\varphi_{d-1}\cdots\varphi_{d-r+1}}
 \cdot \eta^*_i(A^*)\xi_d &= X_i\xi_d,     \\
 \frac{\eta_d(\th_0)}{\vphi}
 \frac{\b{\xi_0,\xi^*_d}}
      {\b{\xi_d,\xi^*_d}}
 \sum_{r=0}^d X_{d-r}E_0\tau^*_{d-r}(A^*)
 \cdot \eta^*_{d-i}(A^*)\xi_d &= X_i\xi_0,    \\
 \frac{1}{\tr(E_dE^*_d)}
 \sum_{r=0}^d 
 \frac{X_{d-r}E_dE^*_d\eta_{r}(A)}
      {\varphi_d\varphi_{d-1}\cdots\varphi_{d-r+1}}
 \cdot \eta^*_{d-i}(A^*)\xi_d &= X_i\xi_d.   
\end{align*}
\end{itemize}
\end{theorem}

\begin{proof}
Apply Theorem \ref{thm:s} to $\Phi^{\downarrow\Downarrow}$.
\end{proof}

\bigskip

{\small

\bibliographystyle{plain}

\begin{thebibliography}{1}

\bibitem{AC}
H.~Alnajjar, B.~Curtin,
A family of tridiagonal pairs,
Linear Algebra Appl.\ 390 (2004) 369--384.

\bibitem{AC2}
H.~Alnajjar, B.~Curtin,
A family of tridiagonal pairs related to the quantum affine 
algebra $U\sb q(\widehat{\mathfrak{sl}}\sb 2)$,
Electron. J. Linear Algebra 13 (2005) 1--9. 

\bibitem{AC3}
H.~Alnajjar, B.~Curtin,
A bilinear form for tridiagonal pairs of $q$-Serre type,
Linear Algebra Appl., 
submitted for publication.

\bibitem{Bas}
P.\ Baseilhac,
A family of tridiagonal pairs and related symmetric functions,
J. Phys. A  39 (38) (2006) 11773--11791.

\bibitem{BT:Borel}
G.~Benkart, P.~Terwilliger,
Irreducible modules for the quantum affine algebra
$U_q(\widehat{\mathfrak{sl}}_2)$ and its Borel subalgebra,
J. Algebra 282 (2004) 172--194;
{\tt arXiv:math.QA/0311152}.

\bibitem{BT:loop}
G.~Benkart, P.~Terwilliger,
The universal central extension of the three-point
$\mathfrak{sl}_2$ loop algebra,
Proc.\ Amer.\ Math.\ Soc.\ 135 (6) (2007) 1659--1668;
{\tt arXiv:math.RA/0512422}.

\bibitem{Bow}
J.\ Bowman,
Irreducible modules for the quantum affine algebra
$U_q(\mathfrak{g})$ and its Borel subalgebra $U_q(\mathfrak{g})^{\geq 0}$,
preprint;
{\tt arxiv:math.QA/0606627}.

\bibitem{C:mlt}
B.\ Curtin,
Modular Leonard triples,
Linear Algebra Appl.,
in press.

\bibitem{C:spinLP}
B.\ Curtin,
Spin Leonard pairs, 
Ramanujan J.\ 13 (2007) 319--332.


\bibitem{F:RL}
D.\ Funk-Neubauer,
Raising/lowering maps and modules for the
quantum affine algebra $U_q(\widehat{\mathfrak{sl}}_2)$,
Communications in Algebra, in press;
{\tt arXiv:math.QA/0506629}.


\bibitem{H}
B.~Hartwig,
Three mutually adjacent Leonard pairs,
Linear Algebra Appl. 408 (2005) 19--39;
{\tt arXiv:math.AC/0508415}.

\bibitem{H:tetra}
B.~Hartwig,
The Tetrahedron algebra and its finite dimensional irreducible modules,
Linear Algebra Appl. 422 (2007) 219--235;
{\tt arXiv:math.RT/0606197}.

\bibitem{HT:tetra}
B.~Hartwig, P.~Terwilliger,
The tetrahedron algebra, the Onsager algebra, 
and the $\mathfrak{sl}_2$ loop algebra,
J. Algebra 308 (2007) 840--863;
{\tt arXiv:math-ph/0511004}.


\bibitem{ITT}
T.~Ito, K.~Tanabe, P.~Terwilliger,
Some algebra related to ${P}$- and ${Q}$-polynomial association schemes, 
Codes and Association Schemes (Piscataway NJ, 1999), 
American Mathematical Society, Providence, RI, 2001, pp. 167--192;
{\tt arXiv:math.CO/0406556}.

\bibitem{IT:shape}
T.~Ito, P.~Terwilliger,
The shape of a tridiagonal pair,
J. Pure Appl.\ Algebra 188 (2004) 145--160;
{\tt arXiv:math.QA/0304244}.

\bibitem{IT:uqsl2hat}
T.~Ito, P.~Terwilliger,
Tridiagonal pairs and the quantum affine algebra
$U_q(\widehat{\mathfrak{sl}}_2)$,
Ramanujan J. 13 (2007) 39--62;
{\tt arXiv:math.QA/0310042}.

\bibitem{IT:non-nilpotent}
T.~Ito, P.~Terwilliger,
\newblock
Two non-nilpotent linear transformations that satisfy the cubic $q$-Serre relations,
J. Algebra Appl., in press;
{\tt arXiv:math.QA/0508398}.

\bibitem{IT:tetra}
T.~Ito, P.~Terwilliger,
\newblock
The $q$-tetrahedron algebra and its finite dimensional irreducible modules,
Comm. Algebra, in press;
{\tt arXiv:math.QA/0602199}.

\bibitem{IT:inverting}
T.~Ito, P.~Terwilliger,
$q$-inverting pairs of linear transformations and the $q$-tetrahedron algebra,
Linear Algebra Appl., submitted for publication;
{\tt arXiv:math.RT/0606237}.

\bibitem{IT:drg}
T.\ Ito, P.\ Terwilliger,
Distance-regular graphs and the $q$-tetrahedron algebra,
European J. Combin.,
submitted for publication;
{\tt arxiv:math.CO/0608694}.

\bibitem{IT:loop}
T.\ Ito, P.\ Terwilliger,
Finite-dimensional irreducible
modules for the three-point $\mathfrak{sl}_2$ loop algebra,
preprint.

\bibitem{IT:Krawt}
T.\ Ito, P.\ Terwilliger,
Tridiagonal pairs of Krawtchouk type,  preprint.

\bibitem{ITW:equitable}
T.~Ito, P.~Terwilliger, C.~Weng,
The quantum algebra $U_q(\mathfrak{sl}_2)$ and its equitable presentation,
J. Algebra 298 (2006) 284--301;
{\tt arXiv:math.QA/0507477}.

\bibitem{M:LT}
S.\ Miklavic,
Leonard triples and hypercubes,
preprint.


\bibitem{N:aw}
K.~Nomura,
\newblock
Tridiagonal pairs and the {A}skey-{W}ilson relations,
Linear Algebra Appl.\ 397 (2005) 99--106.

\bibitem{N:refine}
K.~Nomura,
A refinement of the split decomposition of a tridiagonal pair,
Linear Algebra Appl.\ 403 (2005) 1--23.

\bibitem{N:height1}
K.~Nomura,
Tridiagonal pairs of height one,
Linear Algebra Appl.\ 403 (2005) 118--142.

\bibitem{NT:balanced}
K.~Nomura, P.~Terwilliger,
Balanced Leonard pairs,
Linear Algebra Appl. 420 (2007) 51--69.
{\tt arXiv:math.RA/0506219}. 

\bibitem{NT:formula}
K.~Nomura, P.~Terwilliger,
Some trace formulae involving the split sequences of a Leonard pair,
Linear Algebra Appl.\ 413 (2006) 189--201;
{\tt arXiv:math.RA/0508407}.

\bibitem{NT:det}
K.~Nomura, P.~Terwilliger,
The determinant of $AA^*-A^*A$ for a Leonard pair $A,A^*$,
Linear Algebra Appl.\ 416 (2006) 880--889;
{\tt arXiv:math.RA/0511641}.

\bibitem{NT:mu}
K.~Nomura, P.~Terwilliger,
Matrix units associated with the split basis of a Leonard pair,
Linear Algebra Appl.\ 418 (2006) 775--787;
{\tt arXiv:math.RA/0602416}.

\bibitem{NT:span}
K.~Nomura, P.~Terwilliger,
Linear transformations that are tridiagonal with respect to
both eigenbases of a Leonard pair,
Linear Algebra Appl. 420 (2007) 198--207;
{\tt arXiv:math.RA/0605316}.

\bibitem{NT:switch}
K.~Nomura, P.~Terwilliger,
The switching element for a Leonard pair,
Linear Algebra Appl., submitted for publication;
{\tt arXiv:math.RA/0608623}.

\bibitem{NT:affine}
K.~Nomura, P.~Terwilliger,
Affine transformations of a Leonard pair,
Linear Algebra Appl., submitted for publication;
{\tt arXiv:math.RA/0611783}.

\bibitem{NT:tde}
K.~Nomura, P.~Terwilliger,
The split decomposition of a tridiagonal pair,
Linear Algebra Appl. 424 (2007) 339--345;
{\tt arXiv:math.RA/0612460}.

\bibitem{P}
A.~A.~Pascasio,
On the multiplicities of the primitive idempotents of a
{$Q$}-polynomial distance-regular graph,
European J. Combin.\ 23 (2002) 1073--1078.

\bibitem{R:multi}
H.~Rosengren,
Multivariable orthogonal polynomials and coupling coefficients for 
discrete series representations,
SIAM J. Math. Anal.\ 30 (1999) 233--272.

\bibitem{R:6j}
H.~Rosengren,
An elementary approach to the $6j$-symbols
(classical, quantum, rational, trigonometric, and elliptic),
Ramanujan J. 13 (2007) 131--166;
{\tt arXiv:math.CA/0312310}.


\bibitem{T:sub1}
P.~Terwilliger,
The subconstituent algebra of an association scheme I,
J. Algebraic Combin.\ 1 (1992) 363--388.

\bibitem{T:sub3}
P.~Terwilliger,
The subconstituent algebra of an association scheme III,
J. Algebraic Combin.\ 2 (1993) 177--210.

\bibitem{T:Leonard}
P.~Terwilliger,
Two linear transformations each tridiagonal with respect to an 
eigenbasis of the other,
Linear Algebra Appl.\ 330 (2001) 149--203;
{\tt arXiv:math.RA/0406555}.

\bibitem{T:qSerre}
P.~Terwilliger,
Two relations that generalize the $q$-Serre relations and the
Dolan-Grady relations,
Physics and Combinatorics 1999 (Nagoya), World Scientific Publishing,
River Edge, NJ, 2001, pp.\ 377--398;
{\tt arXiv:math.QA/0307016}.

\bibitem{T:24points}
P.~Terwilliger,
Leonard pairs from 24 points of view,
Rocky Mountain J. Math.\ 32 (2) (2002) 827--888;
{\tt arXiv:math.RA/0406577}.

\bibitem{T:canform}
P.~Terwilliger,
Two linear transformations each tridiagonal with respect to an
eigenbasis of the other; the $TD$-$D$ and the $LB$-$UB$ canonical form,
J. Algebra 291 (2005) 1--45;
{\tt arXiv:math.RA/0304077}.

\bibitem{T:intro}
P.\ Terwilliger,
Introduction to Leonard pairs,
J.\ Comput.\ Appl.\ Math.\ 153 (2) (2003) 463--475.

\bibitem{T:intro2}
P.\ Terwilliger,
Introduction to {L}eonard pairs and {L}eonard systems,
S\=uri\-kaiseki\-kenky\=usho K\=oky\=uroku 1109 (1999) 67--79, 
Algebraic Combinatorics (Kyoto, 1999).

\bibitem{T:split}
P.~Terwilliger,
Two linear transformations each tridiagonal with respect to an
eigenbasis of the other; comments on the split decomposition,
J. Comput.\ Appl.\ Math.\ 178 (2005) 437--452;
{\tt arXiv:math.RA/0306290}.

\bibitem{T:array}
P.~Terwilliger,
Two linear transformations each tridiagonal with respect to an
eigenbasis of the other; comments on the parameter array,
Des.\ Codes Cryptogr.\ 34 (2005) 307--332;
{\tt arXiv:math.RA/0306291}.

\bibitem{T:qRacah}
P.~Terwilliger,
Leonard pairs and the $q$-Racah polynomials,
Linear Algebra Appl.\ 387 (2004) 235--276;
{\tt arXiv:math.QA/0306301}.

\bibitem{T:survey}
P.~Terwilliger,
An algebraic approach to the Askey scheme of orthogonal polynomials,
in: Orthogonal Polynomials and Special Functions: computation and applications,
Lecture Notes in Mathematics, vol.\ 1883, Springer, 2006, pp.\ 255--330;
{\tt arXiv:math.QA/0408390}.

\bibitem{T:Kac-Moody}
P.~Terwilliger,
The equitable presentation for the quantum group $U_q(\mathfrak{g})$ 
associated with a symmetrizable Kac-Moody algebra $\mathfrak{g}$,
J. Algebra 298 (2006) 302--319; 
{\tt arXiv:math.QA/0507478}.

\bibitem{TV}
P.~Terwilliger, R.~Vidunas,
Leonard pairs and the Askey-Wilson relations,
J. Algebra Appl.\ 3 (2004) 411--426;
{\tt arXiv:math.QA/0305356}.

\bibitem{Vidar}
M.\ Vidar,
Tridiagonal pairs of shape $(1,2,1)$, preprint.

\bibitem{V}
R.~Vidunas,
Normalized Leonard pairs and Askey-Wilson relations,
Technical Report MHF 2005-16, Kyushu University 2005;
{\tt arXiv:math.RA/0505041}.

\bibitem{V:AW}
R.~Vidunas,
Askey-Wilson relations and Leonard pairs,
Technical Report MHF 2005-17, Kyushu University 2005;
{\tt arXiv:math.QA/0511509}.

\bibitem{Z}
A.\ Zhedanov,
``Hidden symmetry'' of Askey-Wilson polynomials,
Teoret.\ Mat Fiz.\ 89 (1991) 190--204. 

\bibitem{ZK}
A.\ Zhedanov, A.\ Korovnichenko,
``Leonard pairs'' in classical mechanics,
J. Phys. A  35 (27) (2002) 5767--5780.

\end{thebibliography}

}

\bigskip\bigskip\noindent
Kazumasa Nomura\\
College of Liberal Arts and Sciences\\
Tokyo Medical and Dental University\\
Kohnodai, Ichikawa, 272-0827 Japan\\
email: knomura@pop11.odn.ne.jp

\bigskip\noindent
Paul Terwilliger\\
Department of Mathematics\\
University of Wisconsin\\
480 Lincoln Drive\\ 
Madison, Wisconsin, 53706 USA\\
email: terwilli@math.wisc.edu

\bigskip\noindent
{\bf Keywords.}
Leonard pair, tridiagonal pair, $q$-Racah polynomial, orthogonal polynomial.

\noindent
{\bf 2000 Mathematics Subject Classification}.
05E35, 05E30, 33C45, 33D45.

\end{document}